Konrad-Zuse-Zentrum für Informationstechnik Berlin
Takustr. 7, D-14195 Berlin - Dahlem

Wolfram Koepf

Dieter Schmersau*

# Representations of Orthogonal Polynomials

* Fachbereich Mathematik und Informatik der Freien Universität Berlin

# Representations of Orthogonal Polynomials

Wolfram Koepf
Dieter Schmersau
koepf@zib.de

**Abstract:**

Zeilberger's algorithm provides a method to compute recurrence and differential equations from given hypergeometric series representations, and an adaption of Almquist and Zeilberger computes recurrence and differential equations for hyperexponential integrals. Further versions of this algorithm allow the computation of recurrence and differential equations from Rodrigues type formulas and from generating functions. In particular, these algorithms can be used to compute the differential/difference and recurrence equations for the classical continuous and discrete orthogonal polynomials from their hypergeometric representations, and from their Rodrigues representations and generating functions.

In recent work, we used an explicit formula for the recurrence equation of families of classical continuous and discrete orthogonal polynomials, in terms of the coefficients of their differential/difference equations, to give an algorithm to identify the polynomial system from a given recurrence equation.

In this article we extend these results be presenting a collection of algorithms with which any of the conversions between the differential/difference equation, the hypergeometric representation, and the recurrence equation is possible.

The main technique is again to use explicit formulas for structural identities of the given polynomial systems.

## 1 Structural Formulas for Classical Orthogonal Polynomials

A family

$$y(x) = p_n(x) = k_n\, x^n + \ldots \qquad (n \in \mathbb{N}_0 := \{0, 1, 2, \ldots\}, k_n \neq 0) \tag{1}$$

of polynomials of degree exactly $n$ is a family of classical continuous orthogonal polynomials if it is the solution of a differential equation of the type

$$\sigma(x)\, y''(x) + \tau(x)\, y'(x) + \lambda_n\, y(x) = 0 \tag{2}$$

where $\sigma(x) = ax^2 + bx + c$ is a polynomial of at most second order and $\tau(x) = dx + e$ is a polynomial of first order. Since one demands that $p_n(x)$ has exact degree $n$, by equating the highest coefficients of $x^n$ in (2) one gets

$$\lambda_n = -(an(n-1) + dn)\ . \tag{3}$$

Similarly a family $p_n(x)$ of polynomials of degree exactly $n$, given by (1), is a family of discrete classical orthogonal polynomials if it is the solution of a difference equation of the type

$$\sigma(x)\, \Delta\nabla y(x) + \tau(x)\, \Delta y(x) + \lambda_n\, y(x) = 0\ , \tag{4}$$

where

$$\Delta y(x) = y(x+1) - y(x) \qquad \text{and} \qquad \nabla y(x) = y(x) - y(x-1)$$

denote the forward and backward difference operators, respectively, and $\sigma(x) = ax^2 + bx + c$ and $\tau(x) = dx + e$ are again polynomials of at most second and of first order, respectively. Again, (3) follows.

Since $\Delta\nabla = \Delta - \nabla$, (4) can also be written in the equivalent form

$$\Big(\sigma(x) + \tau(x)\Big)\Delta y(x) - \sigma(x)\,\nabla y(x) + \lambda_n\, y(x) = 0\ ,$$

and replacing $x$ by $x+1$ we arrive at

$$\Big(\sigma(x+1) + \tau(x+1)\Big)\Delta^2 y(x) + \tau(x+1)\,\Delta y(x) + \lambda_n\, y(x+1) = 0\ . \tag{5}$$

It can be shown (see e.g. [14] or [16]) that any solution $p_n(x)$ of either (2) or (4) satisfies a recurrence equation

$$p_{n+1}(x) = (A_n\, x + B_n)\, p_n(x) - C_n\, p_{n-1}(x) \qquad (n \in \mathbb{N}_0\,, p_{-1} \equiv 0) \tag{6}$$

or equivalently

$$x\, p_n(x) = a_n\, p_{n+1}(x) + b_n\, p_n(x) + c_n\, p_{n-1}(x) \tag{7}$$

with

$$a_n = \frac{1}{A_n}\ , \qquad b_n = -\frac{B_n}{A_n}\ , \qquad c_n = \frac{C_n}{A_n}\ . \tag{8}$$

In [14] (compare [16], Eqs. (5) and (10)) we showed that the coefficients $A_n, B_n$, and $C_n$ are given by the explicit formulas

$$A_n = \frac{k_{n+1}}{k_n}\ ,$$

$$B_n = \frac{2\,b\,n\,(\,a\,n + d - a\,) - e\,(\,-d + 2\,a\,)}{(\,d + 2\,a\,n\,)\,(\,d - 2\,a + 2\,a\,n\,)} \cdot \frac{k_{n+1}}{k_n}$$

and

$$C_n = -\frac{\Big((an{+}d{-}2a)n(4ca{-}b^2){+}4a^2c{-}ab^2{+}ae^2{-}4acd{+}db^2{-}bed{+}d^2c\Big)\,(an{+}d{-}2a)n}{(d{-}2a{+}2an)^2\,(2an{-}3a{+}d)\,(2an{-}a{+}d)} \cdot \frac{k_{n+1}}{k_{n-1}}$$

in the continuous case, and by the formulas

$$A_n = \frac{k_{n+1}}{k_n}\ ,$$

$$B_n = \frac{n\,(\,d + 2\,b\,)\,(\,d + a\,n - a\,) + e\,(\,d - 2\,a\,)}{(\,2\,a\,n - 2\,a + d\,)\,(\,d + 2\,a\,n\,)} \cdot \frac{k_{n+1}}{k_n}$$

and

$$\begin{aligned} C_n \;=\; & -\Big((\,n{-}1\,)\,(\,d{+}an{-}a\,)\,(\,a\,n\,d{-}d\,b{-}a\,d{+}a^2n^2{-}2\,a^2\,n{+}4\,c\,a{+}a^2{+}2\,ea{-}b^2\,){-}d\,b\,e{+}d^2\,c{+}a\,e^2\Big) \\ & \cdot \frac{(\,a\,n{+}d{-}2\,a\,)\,n}{(\,d{-}a{+}2\,a\,n)\,(\,d{+}2\,a\,n{-}3\,a\,)\,(\,2\,a\,n{-}2\,a{+}d)^2} \cdot \frac{k_{n+1}}{k_{n-1}} \end{aligned}$$

in the discrete case, in terms of the coefficients $a, b, c, d$ and $e$ of the given differential/difference equation.
Orthogonal polynomials satisfy further structure equations. One of those is given by the derivative/difference rules (see e.g. [14])

$$\sigma(x)\, p_n'(x) = \alpha_n\, p_{n+1}(x) + \beta_n\, p_n(x) + \gamma_n\, p_{n-1}(x) \qquad (n \in \mathbb{N} := \{1,2,3,\ldots\}) \;, \tag{9}$$

and

$$\sigma(x)\, \nabla p_n(x) = \alpha_n\, p_{n+1}(x) + \beta_n\, p_n(x) + \gamma_n\, p_{n-1}(x) \qquad (n \in \mathbb{N}) \;, \tag{10}$$

or

$$\Big(\sigma(x) + \tau(x)\Big)\, \Delta p_n(x) = S_n\, p_{n+1}(x) + T_n\, p_n(x) + R_n\, p_{n-1}(x) \qquad (n \in \mathbb{N}) \;, \tag{11}$$

respectively. Here

$$S_n = \alpha_n \;, \quad T_n = \beta_n - \lambda_n \;, \quad R_n = \gamma_n \;. \tag{12}$$

In [14] we showed that the coefficients $\alpha_n, \beta_n$, and $\gamma_n$ are given by the explicit formulas

$$\alpha_n = a\, n\, \frac{k_n}{k_{n+1}} \;,$$

$$\beta_n = -\frac{n\,(\,a\,n + d - a\,)\,(\,2\,e\,a - d\,b\,)}{(\,d + 2\,a\,n\,)\,(\,d - 2\,a + 2\,a\,n\,)}$$

and

$$\gamma_n = \frac{((n-1)(an+d-a)(4ca-b^2)+ae^2+d^2c-bed)(an+d-a)(an+d-2a)n}{(d-2a+2an)^2\,(2an-3a+d)\,(2an-a+d)} \cdot \frac{k_n}{k_{n-1}}$$

in the continuous case, and by the formulas

$$\alpha_n = a\, n\, \frac{k_n}{k_{n+1}} \;,$$

$$\beta_n = -\frac{n\,(\,d + a\,n - a\,)\,(\,2\,a\,n\,d - a\,d - d\,b + 2\,e\,a - 2\,a^2\,n + 2\,a^2\,n^2\,)}{(\,2\,a\,n - 2\,a + d\,)\,(\,d + 2\,a\,n\,)}$$

and

$$\begin{aligned}\gamma_n \;&=\; \Big((n-1)(d+an-a)(and-db-ad+a^2n^2-2a^2n+4ca+a^2+2ea-b^2)-dbe+d^2c+ae^2\Big)\\ &\quad\cdot \frac{(d+an-a)\,(an+d-2a)\,n}{(d-a+2an)\,(d+2an-3a)\,(2an-2a+d)^2} \cdot \frac{k_n}{k_{n-1}}\end{aligned}$$

in the discrete case, respectively.
Now, we develop further structural identities. Taking the derivative in (2), we get

$$\begin{aligned} 0 \;&=\; \sigma(x)\, p_n'''(x) + \Big(\tau(x) + \sigma'(x)\Big)\, p_n''(x) + \Big(\lambda_n + \tau'(x)\Big)\, p_n'(x) \\ &=\; (ax^2 + bx + c)\, p_n'''(x) + \Big((d + 2a)\, x + (e + b)\Big)\, p_n''(x) + (\lambda_n + d)\, p_n'(x) \;, \end{aligned}$$

hence $y(x) := p_n'(x)$ satisfies a differential equation

$$(a'x^2 + b'x + c')\, y''(x) + (d'x + e')\, y'(x) + \lambda_n'\, y(x)$$

of the same type with

$$a' = a\ , \quad b' = b\ , \quad c' = c\ , \quad d' = d + 2a\ , \quad \text{and} \quad e' = e + b\ . \tag{13}$$

From this we deduce that the equation

$$x\, p_n'(x) = \alpha_n^*\, p_{n+1}'(x) + \beta_n^*\, p_n'(x) + \gamma_n^*\, p_{n-1}'(x)\ , \tag{14}$$

namely a recurrence equation for $p_n'(x)$, is valid, and from (13) it follows that

$$\alpha_n^* = a_n(a, b, c, d + 2a, e + b)\ ,\ \beta_n^* = b_n(a, b, c, d + 2a, e + b)\ ,$$

and

$$\gamma_n^* = c_n(a, b, c, d + 2a, e + b)\ ,$$

where $a_n(a, b, c, d, e), b_n(a, b, c, d, e)$, and $c_n(a, b, c, d, e)$, are given by (8) and the explicit formulas for $A_n, B_n$ and $C_n$.
Similarly in the discrete case, applying $\Delta$ to (4), we get for $y(x) := \Delta p_n(x)$

$$\begin{aligned} 0 &= \Big(\sigma(x+1) - \Delta\sigma(x)\Big)\,\Delta\nabla y(x) + \Big(\tau(x+1) + \Delta\sigma(x)\Big)\,\Delta y(x) + \Big(\lambda_n + \Delta\tau(x)\Big)\,y(x) \\ &= (ax^2 + bx + c)\,\Delta\nabla y(x) + \Big((d + 2a)\,x + d + e + a + b\Big)\,\Delta y(x) + (\lambda_n + d)\,y(x)\ , \end{aligned}$$

hence $y(x) := \Delta p_n(x)$ satisfies a difference equation

$$(a'x^2 + b'x + c')\,\Delta\nabla y(x) + (d'x + e')\,\Delta y(x) + \lambda_n'\, y(x)$$

of the same type with

$$a' = a\ , \quad b' = b\ , \quad c' = c\ , \quad d' = d + 2a\ , \quad \text{and} \quad e' = d + e + a + b\ . \tag{15}$$

From this we deduce that the equation

$$x\,\Delta p_n(x) = \alpha_n^*\,\Delta p_{n+1}(x) + \beta_n^*\,\Delta p_n(x) + \gamma_n^*\,\Delta p_{n-1}(x)\ , \tag{16}$$

namely a recurrence equation for $\Delta p_n(x)$, is valid, and from (15) it follows that

$$\alpha_n^* = a_n(a, b, c, d + 2a, d + e + a + b)\ ,\ \beta_n^* = b_n(a, b, c, d + 2a, d + e + a + b)\ ,$$

and

$$\gamma_n^* = c_n(a, b, c, d + 2a, d + e + a + b)\ ,$$

where $a_n(a, b, c, d, e), b_n(a, b, c, d, e)$, and $c_n(a, b, c, d, e)$, are given by (8) and the explicit formulas for $A_n, B_n$ and $C_n$.
To obtain a derivative rule for $y(x) := p_n'(x)$, we take the derivative of (9) to get

$$\sigma(x)\, p_n''(x) + \sigma'(x)\, p_n'(x) = \alpha_n\, p_{n+1}'(x) + \beta_n\, p_n'(x) + \gamma_n\, p_{n-1}'(x)\ .$$

Applying (14) to replace $x\, p_n'(x)$ results in a derivative rule of the form

$$\sigma(x)\, p_n''(x) = a_n'\, p_{n+1}'(x) + b_n'\, p_n'(x) + c_n'\, p_{n-1}'(x)\ . \tag{17}$$

Similarly in the discrete case a difference rule of the form

$$\sigma(x)\,\Delta\nabla p_n(x) = a'_n\,\Delta p_{n+1}(x) + b'_n\,\Delta p_n(x) + c'_n\,\Delta p_{n-1}(x) \tag{18}$$

can be obtained for $y(x) := \Delta y(x)$.
Finally we substitute (17) in the differential equation. This gives

$$a'_n\,p'_{n+1}(x) + b'_n\,p'_n(x) + c'_n\,p'_{n-1}(x) + \tau(x)\,p'_n(x) + \lambda_n\,p_n(x) = 0\;,$$

and replacing $x\,p'_n(x)$ by (14), again, we obtain an equation of the type

$$p_n(x) = \widehat{a}_n\,p'_{n+1}(x) + \widehat{b}_n\,p'_n(x) + \widehat{c}_n\,p'_{n-1}(x)\;, \tag{19}$$

in the continuous case, and a similar procedure gives

$$p_n(x) = \widehat{a}_n\,\Delta p_{n+1}(x) + \widehat{b}_n\,\Delta p_n(x) + \widehat{c}_n\,\Delta p_{n-1}(x) \tag{20}$$

in the discrete case. Note that in the discrete case also corresponding equations concerning $\nabla$ are valid.
We note in passing that our development shows by simple algebraic arguments that whenever $p_n(x)$ is a polynomial system of degree exactly $n$, satisfying a differential/difference equation of type (2)/(4), a recurrence equation of type (6) and a derivative/difference rule of type (9)/(10) then the system $p'_{n+1}(x)$ $(\Delta p_{n+1}(x))$ is again such a system. This has nothing to do with orthogonality. Indeed, in our further development it will become rather important that in the continuous case the powers $x^n$ and in the discrete case the *falling factorials*

$$x^{\underline{n}} := x\,(x-1)\cdots(x-n+1) = (x-n+1)_n = (-1)^n\,(-x)_n$$

which by no means become orthogonal families, have these properties.
To deduce the coefficients $\alpha^*_n, \beta^*_n, \gamma^*_n$, $a'_n, b'_n, c'_n$, and $\widehat{a}_n, \widehat{b}_n, \widehat{c}_n$, we can follow the above instructions, or we apply the following method: Substituting

$$p_n(x) = k_n\,x^n + k'_n\,x^{n-1} + k''_n\,x^{n-2} + \ldots$$

in the differential/difference equation and equating the coefficients of $x^n$ determines $\lambda_n$, while equating the coefficients of $x^{n-1}$ and $x^{n-2}$ gives $k'_n$, and $k''_n$, respectively, in terms of $k_n$. These values can be substituted in $p_n(x)$. Next, we substitute $p_n(x)$ in the proposed equation, and equate again the three highest coefficients successively to get the three unknowns in terms of $a, b, c, d, e, n, k_{n-1}, k_n$, and $k_{n+1}$ by linear algebra.
These computations can be easily carried out by a computer algebra system, e.g. by Maple. With few seconds of computation time, we get

**Theorem 1** For the solutions of (2) and (4), the relations (14), (17), (19), and (16), (18), (20), respectively, are valid. The coefficients $\alpha^*_n, \beta^*_n, \gamma^*_n$, $a'_n, b'_n, c'_n$, and $\widehat{a}_n, \widehat{b}_n, \widehat{c}_n$, are given by

$$\alpha^*_n = \frac{n}{n+1}\cdot\frac{k_n}{k_{n+1}}\;,$$

$$\beta^*_n = \frac{-2\,b\,n\,(\,a\,n + d - a\,) + d\,(\,b - e\,)}{(\,d + 2\,a\,n\,)\,(\,d - 2\,a + 2\,a\,n\,)}\;,$$

$$\gamma_n^* = -\frac{(\,(\,n-1\,)\,(\,a\,n+d-a\,)\,(\,4\,c\,a-b^2\,)+a\,e^2+d^2\,c-b\,e\,d\,)\,n\,(\,a\,n+d-a\,)}{(\,d-2\,a+2\,a\,n\,)^2\,(\,2\,a\,n-3\,a+d\,)\,(\,2\,a\,n-a+d\,)}\cdot\frac{k_n}{k_{n-1}}\,,$$

$$a_n' = \frac{a\,n\,(n-1)}{n+1}\cdot\frac{k_n}{k_{n+1}}\,,$$

$$b_n' = -\frac{(\,n-1\,)\,(\,a\,n+d\,)\,(\,2\,e\,a-d\,b\,)}{(\,d+2\,a\,n\,)\,(\,d-2\,a+2\,a\,n\,)}\,,$$

$$c_n' = \frac{(\,(\,n-1\,)\,(\,a\,n+d-a\,)\,(\,4\,c\,a-b^2\,)+a\,e^2+d^2\,c-b\,e\,d\,)\,(\,a\,n+d\,)\,(\,a\,n+d-a\,)\,n}{(\,d-2\,a+2\,a\,n\,)^2\,(\,2\,a\,n-3\,a+d\,)\,(\,2\,a\,n-a+d\,)}\cdot\frac{k_n}{k_{n-1}},$$

$$\widehat{a}_n = \frac{1}{n+1}\cdot\frac{k_n}{k_{n+1}}\,,$$

$$\widehat{b}_n = \frac{2\,e\,a-d\,b}{(\,d+2\,a\,n\,)\,(\,d-2\,a+2\,a\,n\,)}\,,$$

$$\widehat{c}_n = \frac{(\,(\,n-1\,)\,(\,a\,n+d-a\,)\,(\,4\,c\,a-b^2\,)+a\,e^2+d^2\,c-b\,e\,d\,)\,a\,n}{(\,d-2\,a+2\,a\,n\,)^2\,(\,2\,a\,n-3\,a+d\,)\,(\,2\,a\,n-a+d\,)}\cdot\frac{k_n}{k_{n-1}}$$

in the continuous case, and

$$\alpha_n^* = \frac{n}{n+1}\cdot\frac{k_n}{k_{n+1}}\,,$$

$$\beta_n^* = \frac{-n\,(\,d+2\,a+2\,b\,)\,(\,d+a\,n-a\,)-d\,(\,e-a-b\,)}{(\,2\,a\,n-2\,a+d\,)\,(\,d+2\,a\,n\,)}\,,$$

$$\begin{aligned}\gamma_n^* \;&=\; -\Big((n-1)\,(d+an-a)\,(and-db-ad+a^2n^2-2a^2n+4ca+a^2+2ea-b^2)-dbe+d^2c+ae^2\Big)\cdot\\ &\quad\frac{(d+an-a)\,n}{(d-a+2an)\,(d+2an-3a)\,(2an-2a+d)^2}\cdot\frac{k_n}{k_{n-1}}\,,\end{aligned}$$

$$a_n' = \frac{a\,n\,(n-1)}{n+1}\cdot\frac{k_n}{k_{n+1}}\,,$$

$$b_n' = -\frac{(\,n-1\,)\,(\,a\,n+d\,)\,(\,2\,a\,n\,d-a\,d-d\,b+2\,e\,a-2\,a^2\,n+2\,a^2\,n^2\,)}{(\,2\,a\,n-2\,a+d\,)\,(\,d+2\,a\,n\,)}\,,$$

$$\begin{aligned}c_n' \;&=\; \Big((n-1)(d+an-a)(and-db-ad+a^2n^2-2a^2n+4ca+a^2+2ea-b^2)-dbe+d^2c+ae^2\Big)\cdot\\ &\quad\frac{(\,d+a\,n-a\,)\,(\,a\,n+d\,)\,n}{(\,d-a+2\,a\,n\,)\,(\,d+2\,a\,n-3\,a\,)\,(\,2\,a\,n-2\,a+d\,)^2}\cdot\frac{k_n}{k_{n-1}},\end{aligned}$$

$$\widehat{a}_n = \frac{1}{n+1}\cdot\frac{k_n}{k_{n+1}}\,,$$

$$\widehat{b}_n = \frac{-2\,a\,n\,(\,d+a\,n-a\,)-d\,b+a\,d-d^2+2\,e\,a}{(\,2\,a\,n-2\,a+d\,)\,(\,d+2\,a\,n\,)}\,,$$

$$\begin{aligned}\widehat{c}_n \;&=\; \Big((n-1)\,(d+an-a)\,(and-db-ad+a^2n^2-2a^2n+4ca+a^2+2ea-b^2)-dbe+d^2c+ae^2\Big)\cdot\\ &\quad\frac{an}{(d-a+2an)\,(d+2an-3a)\,(2an-2a+d)^2}\cdot\frac{k_n}{k_{n-1}}\end{aligned}$$

in the discrete case. □

Note that (19) gives an immediate formula for the antiderivative of a continuous orthogonal polynomial in terms of its neighbors, so that definite integrals can easily be computed, whereas (20) gives an immediate formula for the antidifference of a discrete orthogonal polynomial in terms of its neighbors, so that definite sums can easily be computed.
As a direct consequence of Theorem 1 we have the following representations. The definition of the continuous and discrete families will be given in § 3 and § 5.

**Corollary 1** The classical continuous orthogonal polynomials have the following antiderivative representations:

$$\int H_n(x)\,dx = \frac{1}{2\,(n+1)}\,H_{n+1}(x)$$

(see e.g. [24], (5.5.10)),

$$\int L_n^{(\alpha)}(x)\,dx = L_n^{(\alpha)}(x) - L_{n+1}^{(\alpha)}(x)$$

(see e.g. [25], VI (1.14)),

$$\begin{aligned}\int B_n^{(\alpha)}(x)\,dx &= \frac{2\,(n+1+\alpha)}{(n+1)\,(2n+\alpha+1)\,(2n+\alpha+2)}\,B_{n+1}^{(\alpha)}(x)\\ &+\frac{4}{(2n+\alpha)\,(2n+\alpha+2)}\,B_n^{(\alpha)}(x)\\ &+\frac{2\,n}{(n+\alpha)\,(2n+\alpha)\,(2n+\alpha+1)}\,B_{n-1}^{(\alpha)}(x)\ ,\end{aligned}$$

$$\int C_n^{\alpha}(x)\,dx = \frac{1}{2\,(n+\alpha)}\,C_{n+1}^{\alpha}(x)\,dx - \frac{1}{2\,(n+\alpha)}\,C_{n-1}^{\alpha}(x)\,dx$$

(see e.g. [25], V (7.15)),

$$\begin{aligned}\int P_n^{(\alpha,\beta)}(x)\,dx &= \frac{2\,(n+\alpha+\beta+1)}{(2n+\alpha+\beta+1)\,(2n+\alpha+\beta+2)}\,P_{n+1}^{(\alpha,\beta)}(x)\,dx\\ &+\frac{2\,(\alpha-\beta)}{(2n+\alpha+\beta)\,(2n+\alpha+\beta+2)}\,P_n^{(\alpha,\beta)}(x)\,dx\\ &-\frac{2\,(n+\alpha)\,(n+\beta)}{(n+\alpha+\beta)\,(2n+\alpha+\beta)\,(2n+\alpha+\beta+1)}\,P_{n-1}^{(\alpha,\beta)}(x)\,dx\end{aligned}$$

(see [12], Theorem 6).
The classical discrete orthogonal polynomials have the following antidifference representations:

$$\sum_x c_n^{(\mu)}(x) = -\frac{\mu}{n+1}\,c_{n+1}^{(\mu)}(x)\ ,$$

$$\sum_x k_n^{(p)}(x,N) = k_{n+1}^{(p)}(x,N) - p\,k_n^{(p)}(x,N)\ ,$$

$$\sum_x m_n^{(\gamma,\mu)}(x) = \frac{\mu}{(\mu-1)\,(n+1)}\,m_{n+1}^{(\gamma,\mu)}(x) - \frac{\mu}{\mu-1}\,m_n^{(\gamma,\mu)}(x)\ ,$$

$$\sum_x t_n(x,N) = \frac{1}{2\,(2n+1)}\,t_{n+1}(x,N) - \frac{1}{2}\,t_n(x,N) + \frac{(n-N)\,(n+N)}{2\,(2n+1)}\,t_{n-1}(x,N)\ ,$$

$$\begin{aligned}\sum_x h_n^{(\alpha,\beta)}(x,N) &= \frac{n+\alpha+\beta+1}{(2n+\alpha+\beta+1)\,(2n+\alpha+\beta+2)}\,h_{n+1}^{(\alpha,\beta)}(x,N)\\ &- \frac{2n^2+2n+2n\alpha+2n\beta+\alpha-\alpha N+\beta N+\alpha\beta+\beta+\beta^2}{(2n+\alpha+\beta)\,(2n+\alpha+\beta+2)}\,h_n^{(\alpha,\beta)}(x,N)\\ &+ \frac{(n+\alpha)\,(n+\beta)\,(n-N)\,(n+\alpha+\beta+N)}{(n+\alpha+\beta)\,(2n+\alpha+\beta)\,(2n+\alpha+\beta+1)}\,h_{n-1}^{(\alpha,\beta)}(x,N)\;.\end{aligned}$$

*Proof:* Using the representations for $\hat{a}_n$, $\hat{b}_n$ and $\hat{c}_n$ of Theorem 1 with the particular values for $a, b, c, d, e$ and $k_n$ of the families (see e.g. [1], [19]) give the results. □

Note that the representations for $\hat{a}_n$, $\hat{b}_n$ and $\hat{c}_n$ of Theorem 1, if applied to $P_n(x) = x^n$ or $P_n(x) = x^{\underline{n}}$, respectively, yields the simple results

$$\int x^n\,dx = \frac{1}{n+1} x^{n+1}\;,$$

and

$$\sum_x x^{\underline{n}} = \frac{1}{n+1} x^{\underline{n+1}}\;,$$

respectively. The latter is equivalent to the well-known identity

$$\begin{aligned}\sum_{k=0}^{m} \binom{n+k}{k} &= \frac{1}{n!}\sum_{k=0}^{m}(k+n)^{\underline{n}} = \frac{1}{n!}\sum_{k=n}^{m+n} k^{\underline{n}} = \frac{1}{(n+1)!}\,k^{\underline{n+1}}\Big|_{k=n}^{k=m+n}\\ &= \frac{1}{(n+1)!}\,(n+m)^{\underline{n+1}} = \binom{n+m+1}{m}\;.\end{aligned}$$

The polynomial system

$$K_n^{(\alpha,\beta)}(x) = \left(x+\frac{1+\beta}{\alpha}\right)_n \cdot \alpha^n \cdot {}_1F_1\left(\left.\begin{array}{c} -n \\ 1-x-n-\frac{1+\beta}{\alpha}\end{array}\right| -\frac{1}{\alpha}\right) = (-1)^n \cdot {}_2F_0\left(\left.\begin{array}{c} -n, x+\frac{1+\beta}{\alpha} \\ - \end{array}\right| \alpha\right), \quad (21)$$

which was given in [14], is not orthogonal, but Theorem 1 is still applicable, and we get

$$\sum_x K_n^{(\alpha,\beta)}(x) = \frac{1}{\alpha\,(n+1)}\,K_{n+1}^{(\alpha,\beta)}(x) - K_n^{(\alpha,\beta)}(x)\;.$$

## 2 Connection Coefficients

In this section we would like to consider the problem to determine connection coefficients between different polynomial systems. Here we assume that $P_n(x) = k_n\,x^n + \ldots$ $(n \in \mathbb{N}_0)$ denotes a family of polynomials of degree exactly $n$ and $Q_m(x) = \overline{k}_m\,x^m + \ldots$ $(n \in \mathbb{N}_0)$ denotes a family of polynomials of degree exactly $m$. Then for any $n \in \mathbb{N}_0$ a relation of the type

$$P_n(x) = \sum_{m=0}^{n} C_m(n)\,Q_m(x)\;, \quad (22)$$

is valid, and the coefficients $C_m(n)$ $(n \in \mathbb{N}_0,\ m = 0,\ldots,n)$ are called the *connection coefficients* between the systems $P_n(x)$ and $Q_m(x)$. For simplicity we assume that $C_m(n)$ are defined for all integers $n, m$ and that $C_m(n) = 0$ outside the above $n \times m$-region.

The connection coefficients between many of the classical orthogonal polynomial systems had been determined by different kinds of methods (see e.g. [24], [10], [20]) until Askey and Gasper [6] used recurrence equations to prove the positivity of the connection coefficients between certain instances of the Jacobi polynomials. In a series of papers ([21]–[22], [3]) Ronveaux et. al. recently used such a method more systematically. Here we will present an algorithmic approach different from theirs.

Hence, the main idea is to determine recurrence equations for $C_m(n)$. Since $C_m(n)$ depends on two parameters $m$ and $n$, many mixed recurrence equations are valid as we shall see. The most interesting recurrence equations are those which leave one of the parameters fixed. We will determine those recurrence equations, hence pure recurrence equations w.r.t. $m$ and $n$. The success of this method will heavily depend on whether or not these recurrence equations are of *lowest order*, i.e., whether or not no recurrence equations of lower order for $C_m(n)$ are valid. In cases when the order of the resulting recurrence equation is one, it defines a hypergeometric term which can be given explicitly in terms of *shifted factorials* (or *Pochhammer symbols*) $(a)_k = a(a+1)\cdots(a+k-1) = \Gamma(a+k)/\Gamma(a)$ using the initial value $C_n(n) = k_n/\overline{k}_n$. We will see that there are many instances for this situation.

Note that w.l.o.g. we could assume that $k_n = \overline{k}_m \equiv 1$, i.e., that $\widetilde{P}_n(x)$ and $\widetilde{Q}_m(x)$ are *monic* systems with connection coefficients $\widetilde{C}_m(n)$, because if $P_n(x)$ and $Q_m(x)$ have leading coefficients $k_n$ and $\overline{k}_m$, respectively, then their connection coefficients $C_m(n)$ are given by

$$C_m(n) = \frac{k_n}{\overline{k}_m}\,\widetilde{C}_m(n)\ .$$

In the last section we have already solved a rather special connection problem: (19)/(20) expresses the connection between the polynomial systems $P_n(x) = p_n(x)$ and $Q_m(x) = p'_{m+1}(x)$ or $Q_m(x) = \Delta p_{m+1}(x)$, respectively. In this case the connection coefficients turn out to be rather simple: almost all of them (namely all with $m < n-2$) are zero.

Now, we consider the generic case. We assume that $P_n(x)$ is a polynomial system given by (2)/(4) with $\sigma(x) = ax^2 + bx + c$, and $\tau(x) = dx + e$, and that $Q_m(x)$ is a polynomial system given by (2)/(4) with $\overline{\sigma}(x) = \overline{a}x^2 + \overline{b}x + \overline{c}$, and $\overline{\tau}(x) = \overline{d}x + \overline{e}$.

We know then that both $P_n(x)$ and $Q_m(x)$ satisfy a recurrence equation (7) whose coefficients $a_n(a,b,c,d,e)$, $b_n(a,b,c,d,e)$, and $c_n(a,b,c,d,e)$ were given explicitly in the last section. Note that we will denote all coefficients connected with $Q_m(x)$ by dashes. Hence we have

$$x\,P_n(x) = a_n\,P_{n+1}(x) + b_n\,P_n(x) + c_n\,P_{n-1}(x)$$

and

$$x\,Q_m(x) = \overline{a}_m\,Q_{m+1}(x) + \overline{b}_m\,Q_m(x) + \overline{c}_m\,Q_{m-1}(x)\ ,$$

all of $a_n, b_n, c_n, \overline{a}_m, \overline{b}_m, \overline{c}_m$ given explicitly.

In three steps, we will now derive three independent recurrence equations for $C_m(n)$. First we consider the term $x\,P_n(x)$ (see e.g. [22]). Using the defining equation of $C_m(n)$, and the two recurrence equations for $P_n(x)$ and $Q_m(x)$, we get

$$\begin{aligned}
x\,P_n(x) &= a_n\,P_{n+1}(x) + b_n\,P_n(x) + c_n\,P_{n-1}(x)\\
&= \sum_{m=0}^{n} \Big(a_n\,C_m(n+1)\,Q_m(x) + b_n\,C_m(n)\,Q_m(x) + c_n\,C_m(n-1)\,Q_m(x)\Big)\\
&= \sum_{m=0}^{n} C_m(n)\,x\,Q_m(x)
\end{aligned}$$

$$= \sum_{m=0}^{n} C_m(n) \left( \overline{a}_m \, Q_{m+1}(x) + \overline{b}_m \, Q_m(x) + \overline{c}_m \, Q_{m-1}(x) \right) .$$

By appropriate index shifts, we can equate the coefficient of $Q_m(x)$ to get the "cross rule"

$$a_n \, C_m(n+1) + b_n \, C_m(n) + c_n \, C_m(n-1) = \overline{a}_{m-1} \, C_{m-1}(n) + \overline{b}_m \, C_m(n) + \overline{c}_{m+1} \, C_{m+1}(n) \; . \quad (23)$$

To deduce a second cross rule in terms of the same variables $C_m(n+1), C_m(n), C_m(n-1)$, $C_{m-1}(n)$ and $C_{m+1}(n)$, we examine the term $x \, P_n'(x)$ (or $x \, \Delta P_n(x)$ in the discrete case). Using both recurrence equations for the derivatives/differences

$$x \, P_n'(x) = \alpha_n^* \, P_{n+1}'(x) + \beta_n^* \, P_n'(x) + \gamma_n^* \, P_{n-1}'(x)$$

and

$$x \, Q_m'(x) = \overline{\alpha}_m^* \, Q_{m+1}'(x) + \overline{\beta}_m^* \, Q_m'(x) + \overline{\gamma}_m^* \, Q_{m-1}'(x)$$

(or analogously

$$x \, \Delta P_n(x) = \alpha_n^* \, \Delta P_{n+1}(x) + \beta_n^* \, \Delta P_n(x) + \gamma_n^* \, \Delta P_{n-1}(x)$$

and

$$x \, \Delta Q_m(x) = \overline{\alpha}_m^* \, \Delta Q_{m+1}(x) + \overline{\beta}_m^* \, \Delta Q_m(x) + \overline{\gamma}_m^* \, \Delta Q_{m-1}(x)$$

in the discrete case), we get

$$\begin{aligned} x \, P_n'(x) &= \alpha_n^* \, P_{n+1}'(x) + \beta_n^* \, P_n'(x) + \gamma_n^* \, P_{n-1}'(x) \\ &= \sum_{m=0}^{n} \left( \alpha_n^* \, C_m(n+1) \, Q_m'(x) + \beta_n^* \, C_m(n) \, Q_m'(x) + \gamma_n^* \, C_m(n-1) \, Q_m'(x) \right) \\ &= \sum_{m=0}^{n} C_m(n) \, x \, Q_m'(x) \\ &= \sum_{m=0}^{n} C_m(n) \left( \overline{\alpha}_m^* \, Q_{m+1}'(x) + \overline{\beta}_m^* \, Q_m'(x) + \overline{\gamma}_m^* \, Q_{m-1}'(x) \right) . \end{aligned}$$

Again, by appropriate index shifts, we can equate the coefficient of $Q_m(x)$ to get the cross rule

$$\alpha_n^* \, C_m(n+1) + \beta_n^* \, C_m(n) + \gamma_n^* \, C_m(n-1) = \overline{\alpha}_{m-1}^* \, C_{m-1}(n) + \overline{\beta}_m^* \, C_m(n) + \overline{\gamma}_{m+1}^* \, C_{m+1}(n) \quad (24)$$

(and the same result in the discrete case). In a similar way the cross rule

$$\widehat{a}_n \, C_m(n+1) + \widehat{b}_n \, C_m(n) + \widehat{c}_n \, C_m(n-1) = \overline{\widehat{a}}_{m-1} \, C_{m-1}(n) + \overline{\widehat{b}}_m \, C_m(n) + \overline{\widehat{c}}_{m+1} \, C_{m+1}(n) \quad (25)$$

can be obtained. It turns out, however, that this relation is linearly dependent from (23) and (24), and hence does not yield new information.
Now, we specialize a little. First, we consider the continuous case. To obtain reasonably simple results, we assume furthermore that $\overline{\sigma}(x) = \sigma(x)$. We consider the term $\sigma(x) \, P_n'(x)$. Then, using both derivative rules

$$\sigma(x) \, P_n'(x) = \alpha_n \, P_{n+1}(x) + \beta_n \, P_n(x) + \gamma_n \, P_{n-1}(x)$$

and

$$\overline{\sigma}(x)\, Q'_m(x) = \overline{\alpha}_m\, Q_{m+1}(x) + \overline{\beta}_m\, Q_m(x) + \overline{\gamma}_m\, Q_{m-1}(x)\ ,$$

we get

$$\begin{aligned}
\sigma(x)\, P'_n(x) &= \alpha_n\, P_{n+1}(x) + \beta_n\, P_n(x) + \gamma_n\, P_{n-1}(x) \\
&= \sum_{m=0}^{n} \left(\alpha_n\, C_m(n+1)\, Q_m(x) + \beta_n\, C_m(n)\, Q_m(x) + \gamma_n\, C_m(n-1)\, Q_m(x)\right) \\
&= \sum_{m=0}^{n} C_m(n)\, \sigma(x)\, Q'_m(x) \\
&= \sum_{m=0}^{n} C_m(n) \left(\overline{\alpha}_m\, Q_{m+1}(x) + \overline{\beta}_m\, Q_m(x) + \overline{\gamma}_m\, Q_{m-1}(x)\right) .
\end{aligned}$$

Again, by appropriate index shifts, this results in the cross rule

$$\alpha_n\, C_m(n+1) + \beta_n\, C_m(n) + \gamma_n\, C_m(n-1) = \overline{\alpha}_{m-1}\, C_{m-1}(n) + \overline{\beta}_m\, C_m(n) + \overline{\gamma}_{m+1}\, C_{m+1}(n)\ . \quad (26)$$

To obtain a pure recurrence equation w.r.t. $m$, from the three cross rules (23), (24), and (26) by linear algebra we eliminate the variables $C_m(n+1)$ and $C_m(n-1)$, and to obtain a pure recurrence equation w.r.t. $n$, we eliminate the variables $C_{m-1}(n)$ and $C_{m+1}(n)$. For simplicity we consider the monic case.

**Theorem 2** Let $P_n(x)$ be a monic polynomial system given by the differential equation (2) with $\sigma(x) = ax^2 + bx + c$, and $\tau(x) = dx + e$, and let $Q_m(x)$ be a monic polynomial system given by (2) with $\overline{\sigma}(x) = \sigma(x)$, and $\overline{\tau}(x) = \overline{d}x + \overline{e}$. Then the relation (22) is valid, $C_m(n)$ satisfying the recurrence equation

$$\begin{aligned}
&-(m-n)(am+d-a+an)(\overline{d}+2am)(\overline{d}+a+2am)(\overline{d}+3a+2am) \\
&\quad (\overline{d}+2am+2a)^2\, C_m(n) + (-d\,b\,n\,\overline{d} + 2\,d\,a\,m^2\,b + d\,b\,\overline{d} + 2\,d\,a\,m\,b + 2\,d\,\overline{e}\,n\,a \\
&\quad + d\,\overline{d}\,\overline{e} + 2\,d\,\overline{d}\,b\,m - m\,b\,\overline{d}^2 - e\,\overline{d}^2 - 4\,a^2\,m^2\,e - m^2\,a\,b\,\overline{d} + b\,n\,\overline{d}\,a - 2\,e\,\overline{d}\,a \\
&\quad - 4\,a^2\,m\,e - 4\,e\,\overline{d}\,a\,m + 2m^2\,a^2\,\overline{e} + 2\,\overline{e}\,a^2\,n^2 - 2\,\overline{e}\,a^2\,n - m\,a\,b\,\overline{d} + 2\,m\,\overline{d}\,\overline{e}\,a \\
&\quad + 2\,m\,\overline{e}\,a^2 - b\,n^2\,\overline{d}\,a)(\overline{d}+2am+2a)(m+1)(\overline{d}+a+2am)(\overline{d}+3a+2am) \\
&\quad C_{m+1}(n) - (\overline{d}+2am)(m+1)(-am-2a+an-\overline{d}+d)(am+an+a+\overline{d}) \\
&\quad (a\,b^2\,m^2 - 4\,a^2\,m^2\,c - 8\,a^2\,m\,c + 2\,a\,m\,b^2 - 4\,a\,\overline{d}\,m\,c + m\,b^2\,\overline{d} - 4\,a\,\overline{d}\,c - a\,\overline{e}^2 + a\,b^2 - c\,\overline{d}^2 \\
&\quad + b\,\overline{e}\,\overline{d} - 4\,a^2\,c + b^2\,\overline{d})(m+2)\, C_{m+2}(n) = 0
\end{aligned}$$

with respect to $m$, with initial values $C_n(n) = 1, C_{n+1}(n) \equiv 0$. Furthermore the recurrence equation

$$\begin{aligned}
&-(d+2an)^2(d-a+2an)(d+2an+2a)(d+a+2an)(-m+n+2) \\
&\quad (\overline{d}+am+a+an)\, C_m(n+2) + (d-a+2an)(d+a+2an)(n+2)(d+2an) \\
&\quad (-2\,\overline{e}\,a\,d - b\,d^2\,n - 2\,m\,a^2\,e + 2\,m^2\,a^2\,e + 2\,a^2\,e\,n - 4\,\overline{e}\,a^2\,n^2 + 2\,\overline{d}\,b\,d - 2\,\overline{d}\,e\,a \\
&\quad + m\,b\,d\,a - b\,d^2 + 2\,e\,d\,a - 4\,\overline{e}\,a\,d\,n + 2\,e\,d\,a\,n - a\,b\,d\,n^2 - b\,d\,n\,a + 2\,a^2\,e\,n^2 - m^2\,a\,b\,d \\
&\quad - 4\,\overline{e}\,a^2\,n - \overline{d}\,m\,b\,d + 2\,\overline{d}\,a\,n\,b + 2\,\overline{d}\,a\,n^2\,b + 2\,\overline{d}\,d\,b\,n + 2\,\overline{d}\,m\,e\,a \\
&\quad - \overline{e}\,d^2 + \overline{d}\,d\,e)C_m(n+1) + (d+2an+2a)(n+2)(n+1)(an-am+d-\overline{d}) \\
&\quad (an+am-a+d)(b\,e\,d - a\,e^2 - d^2\,c - 4\,a\,c\,n\,d - 4\,a^2\,c\,n^2 + a\,b^2\,n^2 + n\,b^2\,d)\, C_m(n) = 0
\end{aligned}$$

with respect to $n$ is valid.

*Proof:* Using the explicit representations given in the last section in combination with (23), (24), and (26), and elimination of $C_m(n+1)$ and $C_m(n-1)$, or $C_{m-1}(n)$ and $C_{m+1}(n)$, respectively, yields the results. □

Note that the recurrence equation given in Theorem 2 reduces to two terms, and hence can be represented by hypergeometric terms, for the connection between Laguerre polynomials ($P_n(x) = L_n^{(\alpha)}(x), Q_m(x) = L_m^{(\beta)}(x)$), and between the Gegenbauer polynomials ($P_n(x) = C_n^{\alpha}(x), Q_m(x) = C_m^{\beta}(x)$). We will consider these and more cases by another method in § 4.

Now, let's switch to the discrete case. There are two possibilities to obtain a relation similar to (26). Replacing the derivative by $\nabla$, the same argument gives (26), again, valid for $\sigma(x) = \overline{\sigma}(x)$. If $\sigma(x) + \tau(x) = \overline{\sigma}(x) + \overline{\tau}(x)$, we can replace the derivative by $\Delta$, and adopt the above argument to get the relation

$$S_n\, C_m(n+1) + T_n\, C_m(n) + R_n\, C_m(n-1) = \overline{S}_{m-1}\, C_{m-1}(n) + \overline{T}_m\, C_m(n) + \overline{R}_{m+1}\, C_{m+1}(n)\ . \quad (27)$$

Hence we get

**Theorem 3** Let $P_n(x)$ be a monic polynomial system given by the difference equation (4) with $\sigma(x) = ax^2 + bx + c$, and $\tau(x) = dx + e$, and let $Q_m(x)$ be a monic polynomial system given by (4) with $\overline{\sigma}(x) = \sigma(x)$, and $\overline{\tau}(x) = \overline{d}x + \overline{e}$. Then the relation (22) is valid, $C_m(n)$ satisfying the recurrence equation

$$\begin{aligned}
&(\overline{d} + 2am + 2a)^2 (\overline{d} + 3a + 2am)(\overline{d} + a + 2am)(\overline{d} + 2am)(-m + n)\\
&\quad (an - a + d + am)\, C_m(n) - (\overline{d} + 2am + 2a)(\overline{d} + 3a + 2am)(\overline{d} + a + 2am)(m + 1)\\
&\quad (-a\overline{d}d + an^2\overline{d}^2 + 2ea\overline{d} - 2a^2\overline{e}m^2 + b\overline{d}^2 m - 2na^3m^2 - 2a^3nm - a^2n\overline{d}\\
&\quad + a^2n^2\overline{d} - an\overline{d}^2 - 2\overline{e}a\overline{d}m + 4ea^2m + ba\overline{d}m^2 + 2a^3m^2n^2 + a\overline{d}bm\\
&\quad + 2a^2m\overline{d}n^2 - 4a^3m^3 - 2a\overline{d}^2m^2 - 4a^2\overline{d}m^3 + dn\overline{d}^2 - 2a^2m^2d - a^2m\overline{d} - \overline{d}^2ma\\
&\quad - 2a^2\overline{e}m + 2a^3mn^2 + and\overline{d} - 2a^3m^2 - 2a^2\overline{e}n^2 - \overline{d}^2d - \overline{d}bd - 2a^3m^4 - 2a^2md\\
&\quad + 2\overline{d}mand + 2a^2m^2nd + 4a^2em^2 + 4aem\overline{d} - 2a\overline{e}nd + 2a^2ndm - 2a^2m\overline{d}n\\
&\quad - 3am\overline{d}d - \overline{d}m^2ad - 2\overline{d}bmd - 2am^2bd - 2ambd - 5\overline{d}m^2a^2 + 2\overline{e}a^2n\\
&\quad - \overline{d}\overline{e}d - \overline{d}^2md + e\overline{d}^2 + dbn\overline{d} + an^2b\overline{d} - anb\overline{d})\, C_{m+1}(n) + (m + 1)\\
&\quad (\overline{d} + 2am)(\overline{d} + am + a + an)(m + 2)(4a^2cm^2 + 2a^2\overline{e}m^2 - b\overline{d}^2m - b^2am^2\\
&\quad + 2\overline{e}a\overline{d}m - ba\overline{d}m^2 + 4\overline{d}cam + 2a^2\overline{d} + 4a^3m - 2a\overline{d}bm + 4a^3m^3\\
&\quad - \overline{d}b^2m + a\overline{d}^2m^2 + 2a^2\overline{d}m^3 + 6a^2m\overline{d} + 2\overline{d}^2ma + 4a^2\overline{e}m + 6a^3m^2\\
&\quad + 4\overline{d}ca - 2b^2am + 8a^2cm + 2a^2\overline{e} - b^2a + a^3m^4 + a^3 + 4a^2c - \overline{d}b^2 - ba\overline{d} + \overline{d}^2c\\
&\quad + a\overline{d}^2 - b\overline{d}^2 + a\overline{e}^2 - \overline{d}b\overline{e} + 2\overline{e}a\overline{d} + 6\overline{d}m^2a^2)\\
&\quad (-am - 2a - \overline{d} + an + d)\, C_{m+2}(n) = 0
\end{aligned}$$

with respect to $m$, with initial values $C_n(n) = 1, C_{n+1}(n) \equiv 0$. Furthermore the recurrence equation

$$\begin{aligned}
&(d + 2an + 2a)(n + 2)(-n^2ab^2 - d^2bn - dbe + 2a^2n^2e + 4n^2a^2c\\
&\quad + 2dn^3a^2 + d^2an^2 - db^2n + a^3n^4 + d^2c + ae^2 - dan^2b + 2dane + 4dcna)\\
&\quad (n + 1)(an - a + d + am)(-am - \overline{d} + an + d)\, C_m(n) - (d + 2an)\\
&\quad (d - a + 2an)(d + a + 2an)(n + 2)(-2ea\overline{d} - 2na^3m^2 + ed\overline{d}
\end{aligned}$$

$$+ 2 a^3 n m + a m d^2 - a m^2 d^2 + d^2 n \overline{d} - 2 e a^2 m + a n^2 d \overline{d} - 2 a^3 m^2 n^2$$
$$- 2 a^2 m \overline{d} n^2 - d^2 b n - a^2 m^2 d - a n d b + 2 a^2 n^2 e + 2 a^3 m n^2 + 3 a^2 n d$$
$$+ 7 a^2 n^2 d + a n d \overline{d} + 2 e a d - \overline{d} m d^2 + 4 d n^3 a^2 - 4 a^2 \overline{e} n^2$$
$$- \overline{e} d^2 + 2 a^2 n e + d^2 \overline{d} + 2 d^2 a n^2 + 2 \overline{d} b d + 2 a^3 n^4 + a d^2 - b d^2$$
$$- d a n^2 b + a^2 m d + 2 d a n e + 4 a^3 n^3 + 2 a^3 n^2 - 2 \overline{d} m a n d - 2 a^2 m^2 n d$$
$$+ 2 a^2 e m^2 + 2 a e m \overline{d} - 4 a \overline{e} n d - 2 a \overline{e} d + 2 a^2 n d m$$
$$- 2 a^2 m \overline{d} n - a m \overline{d} d - \overline{d} b m d - a m^2 b d + a m b d - 4 \overline{e} a^2 n$$
$$+ 2 d b n \overline{d} + 2 a n^2 b \overline{d} + 2 a n b \overline{d} + 3 a n d^2 ) C_m(n+1) +$$
$$( d + 2 a n )^2 ( d + 2 a n + 2 a ) ( d - a + 2 a n ) ( d + a + 2 a n ) ( -m + n + 2 )$$
$$( \overline{d} + a m + a + a n ) C_m(n+2) = 0$$

with respect to $n$ is valid.
Next let $\overline{\sigma}(x) + \overline{\tau}(x) = \sigma(x) + \tau(x)$, hence $\overline{a} = a$, $\overline{b} = b + f$, $\overline{c} = c + g$, $\overline{d} = d - f$, $\overline{e} = e - g$ for some constants $f, g$. Then the relation (22) is valid, $C_m(n)$ satisfying the recurrence equation

$$( -d + f - 2 a m ) ( -d + f - 2 a m - 2 a )^2 ( -d + f - a - 2 a m ) ( -d + f - 3 a - 2 a m )$$
$$( -m + n ) ( a n - a + d + a m ) C_m(n) - ( -d + f - 2 a m - 2 a ) ( -d + f - a - 2 a m )$$
$$( -d + f - 3 a - 2 a m ) ( m + 1 )(2 e a^2 m - 2 a^3 m^2 n - d^3 + 2 a^2 g m + 2 e a d - a d^2$$
$$- b d^2 + d^2 b n + d^2 a n^2 + a^2 n^2 d - 2 a^2 n^2 e + 2 a^2 n e - a^2 n d + d a n^2 b - a n d b$$
$$- 2 d a n e - 2 a e f - d a n^2 f + 2 a^3 m^2 n^2 + 2 a n m d^2 - 2 a^3 m^2 - a m^2 b d$$
$$- a m b d - 7 a^2 m^2 d - 3 a^2 m d - 3 a m^2 d^2 - 4 a m d^2 + f^3 m - 4 a^3 m^3 - 2 a^3 m^4$$
$$- a m f b - 2 a^2 m f n^2 + 2 a^2 n^2 d m + 2 a^2 f m n - 2 a^3 m n + 2 a^3 n^2 m - m d^3$$
$$+ f b a n + 2 d g a n + 2 a^2 m^2 n d + d^3 n + 2 a^2 e m^2 + a^2 n f - d b n f - d^2 n f$$
$$- a^2 n^2 f - a n^2 b f - 2 g a^2 n + 2 g a^2 n^2 - 2 m f a n d + m f d^2 + 3 f a d - 2 f^2 a$$
$$- 2 f^2 d + 2 f d^2 + 4 a m^2 f d + 2 d e a m + d^2 g + 5 a^2 m f + 9 a^2 m^2 f - m d^2 b$$
$$- f e d - f g d - m f^2 d + 8 a m f d + f^2 e - 2 a e m f + 2 a d g m - 2 a f g m$$
$$- 4 a^2 d m^3 + 4 a^2 f m^3 + d f b + f^2 b m - 3 a m^2 f^2 + 2 a^2 m^2 g - a m^2 f b - 6 f^2 a m$$
$$+ f^3 ) C_{m+1}(n) - ( -d + f - 2 a m ) ( m + 1 )(4 e a^2 m + 8 a^2 c m - 2 b^2 a m$$
$$+ 4 a^2 g m + 2 e a d + 4 d c a - d b^2 + d^2 c + a d^2 - b d^2 - b^2 a + a e^2 + 2 a^2 e$$
$$+ 2 a^2 d + 4 a^2 c + a^3 - d b e - b a d - 2 a e f + 6 a^3 m^2 + 4 a^3 m - a m^2 b d$$
$$- 2 a m b d + 6 a^2 m^2 d + 6 a^2 m d + a m^2 d^2 + 2 a m d^2 + 4 a^3 m^3 + a^3 m^4$$
$$- b^2 a m^2 + 4 a^2 c m^2 - 2 a m f b + 2 a^2 e m^2 - m f d^2 - 3 f a d + f^2 a + f^2 d - f d^2$$
$$- 3 a m^2 f d + 2 d e a m + d^2 g - 6 a^2 m f - 6 a^2 m^2 f - m d^2 b + 2 d g a - f e d$$
$$- 2 f g a - f g d + m f^2 d - 6 a m f d + f^2 e - 2 a e m f + a g^2 + 4 a d c m + 2 a d g m$$
$$- 4 a f c m - 2 a f g m + f^2 c - 2 d f c - 2 a e g + d b g + f b e - f b g + 2 a^2 d m^3$$
$$- 2 a^2 f m^3 - d b^2 m + f b^2 m + f^2 b m + a m^2 f^2 + 2 a^2 m^2 g + f^2 b - a m^2 f b - a f b$$
$$+ 2 f^2 a m - 2 a^2 f + f b^2 + 2 a^2 g - 4 a f c)( m + 2 ) ( -d + f - a m - a n - a )$$
$$( -a m - 2 a + f + a n ) C_{m+2}(n) = 0$$

with respect to $m$, with initial values $C_n(n) = 1, C_{n+1}(n) \equiv 0$. Furthermore the recurrence equation

$$( n + 2 ) ( d + 2 a n + 2 a )(2 d n^3 a^2 + a^3 n^4 + d^2 c + a e^2 - d b e - d b^2 n - d^2 b n$$

$$- n^2 a b^2 + d^2 a n^2 + 2 a^2 n^2 e + 4 n^2 a^2 c - d a n^2 b + 4 d c n a + 2 d a n e)$$
$$(n+1)(an - a + d + am)(-am + f + an) C_m(n) - (d - a + 2an)(d + a + 2an)$$
$$(d + 2an)(n+2)(-2 e a^2 m - 2 a^3 m^2 n + d^3 + 4 d n^3 a^2 - 2 e a d + 2 a^3 n^4$$
$$+ a d^2 + b d^2 + 4 a^3 n^3 + 2 a^3 n^2 + d^2 b n + 3 d^2 a n^2 + 4 a n d^2 + 7 a^2 n^2 d$$
$$- 2 a^2 n^2 e - 2 a^2 n e + 3 a^2 n d + d a n^2 b + a n d b - 2 d a n e + 2 a e f - d a n^2 f$$
$$- 2 a^3 m^2 n^2 - 2 a n m d^2 - a m^2 b d + a m b d - a^2 m^2 d + a^2 m d - a m^2 d^2$$
$$+ 2 a^2 m f n^2 - 2 a^2 n^2 d m + 2 a^2 f m n + 2 a^3 m n + 2 a^3 n^2 m - m d^3 - 2 f b a n$$
$$- d f a n + 4 d g a n - 2 a^2 m^2 n d + d^3 n + 2 a^2 e m^2 - 2 d b n f - d^2 n f - 2 a n^2 b f$$
$$+ 4 g a^2 n + 4 g a^2 n^2 + 2 m f a n d + m f d^2 + m f b d - f d^2 + 2 d e a m + d^2 g$$
$$- m d^2 b + 2 d g a - f e d + a m f d - 2 a e m f - 2 d f b) C_m(n+1) - (d - a + 2an)$$
$$(d + a + 2an)(d + 2an)^2 (d + 2an + 2a)(-m + n + 2)(-d + f - am - an - a)$$
$$C_m(n+2) = 0$$

with respect to $n$ is valid. □

Note that the recurrence equation for $\overline{\sigma}(x) = \sigma(x)$ given in Theorem 3 reduces to two terms, and hence can be represented by hypergeometric terms, for the connection between Charlier polynomials ($P_n(x) = c_n^{(\mu)}(x), Q_m(x) = c_m^{(\nu)}(x)$), between Meixner polynomials ($P_n(x) = m_n^{(\gamma,\mu)}(x), Q_m(x) = m_m^{(\delta,\mu)}(x)$), and between Krawchouk polynomials ($P_n(x) = k_n^{(p)}(x, N)$, $Q_m(x) = k_m^{(p)}(x, M)$), We will consider these and more cases by another method in § 6.

# 3 Hypergeometric Representations: Continuous Case

Note that by $\widetilde{P}_n^{(\alpha,\beta)}(x), \widetilde{C}_n^{\alpha}(x), \widetilde{L}_n^{(\alpha)}(x), \widetilde{H}_n(x), \widetilde{B}_n^{(\alpha)}(x)$ we denote the monic *Jacobi, Gegenbauer, Laguerre, Hermite* and *Bessel polynomials*. Their non-monic counterparts have the standardizations (see [1], (22.3), and [2]; Al-Salam denotes the Bessel polynomials by $Y_n^{(\alpha)}(x)$)

| system | $P_n^{(\alpha,\beta)}(x)$ | $C_n^{\alpha}(x)$ | $L_n^{(\alpha)}(x)$ | $H_n(x)$ | $B_n^{(\alpha)}(x)$ |
|---|---|---|---|---|---|
| $k_n$ | $\frac{1}{2^n}\binom{2n+\alpha+\beta}{n}$ | $\frac{(\alpha)_n\, 2^n}{n!}$ | $\frac{(-1)^n}{n!}$ | $2^n$ | $\frac{(n+\alpha+1)_n}{2^n}$ |

We get

**Theorem 4** Let $P_n(x)$ be a monic polynomial system given by the differential equation (2) with $\sigma(x) = ax^2 + bx + c$, and $\tau(x) = dx + e$. Then the power series coefficients $C_m(n)$ given by

$$P_n(x) = \sum_{m=0}^{n} C_m(n)\, x^m \tag{28}$$

satisfy the recurrence equation

$$(m-n)(an+d-a+am)C_m(n)+(m+1)(bm+e)C_{m+1}(n)+c(m+1)(m+2)C_{m+2}(n) = 0\ . \tag{29}$$

In particular, if $c = 0$, then the recurrence equation

$$(m-n)(an + d - a + am)C_m(n) + (m+1)(bm + e)C_{m+1}(n) = 0$$

is valid, and we have the hypergeometric representation

$$P_n(x) = \frac{\left(\frac{e}{b}\right)_n \left(\frac{d-a}{a}\right)_n}{\left(\frac{d-a}{2a}\right)_n \left(\frac{d}{2a}\right)_n} \left(\frac{b}{4a}\right)^n \cdot {}_2F_1\left(\begin{array}{c} -n, \frac{d+(n-1)a}{a} \\ \frac{e}{b} \end{array} \middle| -\frac{a}{b}\,x\right), \tag{30}$$

valid for $a \neq 0$, or

$$P_n(x) = \left(\frac{e}{b}\right)_n \left(\frac{b}{d}\right)^n {}_1F_1\left(\begin{array}{c} -n \\ \frac{e}{b} \end{array} \middle| -\frac{d}{b}\,x\right), \tag{31}$$

valid for $a = 0, b \neq 0$, or finally

$$P_n(x) = \left(\frac{e}{d}\right)^n {}_1F_0\left(\begin{array}{c} -n \\ - \end{array} \middle| -\frac{d}{e}\,x\right), \tag{32}$$

valid for $a = 0, b = 0$.

Therefore, the classical continuous orthogonal polynomials and their monic counterparts have the following hypergeometric power series representations:

$$\begin{aligned}
P_n^{(\alpha,\beta)}(x) &= \binom{n+\alpha}{n}\, {}_2F_1\left(\begin{array}{c} -n, n+\alpha+\beta+1 \\ \alpha+1 \end{array} \middle| \frac{1-x}{2}\right) \\
&= \binom{2n+\alpha+\beta}{n} \left(\frac{x-1}{2}\right)^n {}_2F_1\left(\begin{array}{c} -n, -n-\alpha \\ -2n-\alpha-\beta \end{array} \middle| \frac{2}{1-x}\right) \\
&= (-1)^n \binom{n+\beta}{n}\, {}_2F_1\left(\begin{array}{c} -n, n+\alpha+\beta+1 \\ \beta+1 \end{array} \middle| \frac{1+x}{2}\right) \\
&= \binom{2n+\alpha+\beta}{n} \left(\frac{x+1}{2}\right)^n {}_2F_1\left(\begin{array}{c} -n, -n-\beta \\ -2n-\alpha-\beta \end{array} \middle| \frac{2}{1+x}\right),
\end{aligned}$$

$$\widetilde{C}_n^{\alpha}(x) = x^n\, {}_2F_1\left(\begin{array}{c} -n/2, -n/2+1/2 \\ -n-\alpha+1 \end{array} \middle| \frac{1}{x^2}\right),$$

$$C_n^{\alpha}(x) = \frac{(\alpha)_n\, 2^n\, x^n}{n!}\, {}_2F_1\left(\begin{array}{c} -n/2, -n/2+1/2 \\ -n-\alpha+1 \end{array} \middle| \frac{1}{x^2}\right),$$

$$\widetilde{L}_n^{(\alpha)}(x) = (1+\alpha)_n\, (-1)^n\, {}_1F_1\left(\begin{array}{c} -n \\ 1+\alpha \end{array} \middle| x\right) = x^n\, {}_2F_0\left(\begin{array}{c} -n, -n-\alpha \\ - \end{array} \middle| -\frac{1}{x}\right),$$

$$L_n^{(\alpha)}(x) = \binom{n+\alpha}{n}\, {}_1F_1\left(\begin{array}{c} -n \\ 1+\alpha \end{array} \middle| x\right) = \frac{(-x)^n}{n!}\, {}_2F_0\left(\begin{array}{c} -n, -n-\alpha \\ - \end{array} \middle| -\frac{1}{x}\right),$$

$$\widetilde{H}_n(x) = x^n\, {}_2F_0\left(\begin{array}{c} -n/2, -n/2+1/2 \\ - \end{array} \middle| -\frac{1}{x^2}\right),$$

$$H_n(x) = 2^n\, x^n\, {}_2F_0\left(\begin{array}{c} -n/2, -n/2+1/2 \\ - \end{array} \middle| -\frac{1}{x^2}\right),$$

$$\widetilde{B}_n^{(\alpha)}(x) = \frac{2^n}{(n+\alpha+1)_n}\, {}_2F_0\left(\begin{array}{c} -n, n+\alpha+1 \\ - \end{array} \middle| -\frac{x}{2}\right) = x^n\, {}_1F_1\left(\begin{array}{c} -n \\ -2n-\alpha \end{array} \middle| \frac{2}{x}\right).$$

$$B_n^{(\alpha)}(x) = {}_2F_0\left(\begin{array}{c} -n, n+\alpha+1 \\ - \end{array} \middle| -\frac{x}{2}\right) = \frac{(n+\alpha+1)_n}{2^n}\, x^n\, {}_1F_1\left(\begin{array}{c} -n \\ -2n-\alpha \end{array} \middle| \frac{2}{x}\right).$$

These results are all particular cases of the recurrence equation (29).

*Proof:* Substituting the power series (28) into the differential equation, and equating the coefficients yields the recurrence equation (29).
For $c = 0$ this recurrence equation degenerates to a two-term recurrence equation, and hence establishes the hypergeometric representations (30)–(32), using the initial value $C_n(n) = 1$.
A shift in the $x$-variable then generates the representations for the Jacobi polynomials. The two points of development $x_1 = 1$ and $x_2 = -1$ correspond to the zeros of $\sigma(x)$.
Note that some of the hypergeometric representations correspond to each other by changing the direction of summation.
The other representations follow by substituting the particular parameters $a, b, c, d$, and $e$ into the recurrence equation (29), and using the initial value $C_n(n) = k_n$ (or $C_n(n) = 1$ in the monic case). □

We would like to mention that the recurrence equation (29) carries the complete information about the hypergeometric representations given in the theorem.
The method described results in four different hypergeometric representations for the Jacobi polynomials. Many more hypergeometric representations exist, but the algorithmic procedure presented finds *power series* representations only. For example, the representation (see e.g. [1] (22.5.45))

$$P_n^{(\alpha,\beta)}(x) = \binom{n+\beta}{n} \left(\frac{x-1}{2}\right)^n {}_2F_1\left(\left.\begin{array}{c} -n, -n-\alpha \\ \beta+1 \end{array}\right| \frac{x+1}{x-1}\right)$$

cannot be discovered by this method.
The method was able to find hypergeometric series representations with point of development $x_0 = 0$ for the Gegenbauer polynomials which are specific Jacobi polynomials, but failed in the Jacobi case, though. One might ask whether such a representation exists. This question can be completely answered by an algorithm of Petkovšek [18]. Petkovšek's algorithm finds *all* hypergeometric term solutions of holonomic recurrence equations, i.e., homogeneous linear recurrence equations with polynomial coefficients. Using the recurrence equation (29), an application of Petkovšek's algorithm *proves* that the Jacobi polynomials do *not* generally have a hypergeometric series representation at the origin.
Note that the method of the last section, although more complicated, does also give the recurrence equation (29), and hence the above results.

# 4 Power Representations

Whereas in the last section we considered the specific connection coefficient problem for $Q_m(x) = x^m$, in this section the opposite problem, having $P_n(x) = x^n$, is studied.
In many applications, one wants to develop a given polynomial in terms of a given orthogonal polynomial system. In this case handy formulas for the powers $x^n$ are very welcome.

**Theorem 5** Let $Q_m(x)$ be a monic polynomial system given by the differential equation (2) with $\overline{\sigma}(x) = \overline{a}x^2 + \overline{b}x + \overline{c}$, and $\overline{\tau}(x) = \overline{d}x + \overline{e}$. Then the coefficients $C_m(n)$ of the power representations

$$x^n = \sum_{m=0}^{n} C_m(n)\, Q_m(x)$$

satisfy the recurrence equation

$$(n-m)\,(\overline{d}+2\overline{a}m)\,(\overline{d}+3\overline{a}+2\overline{a}m)\,(\overline{d}+\overline{a}+2\overline{a}m)\,(\overline{d}+2\overline{a}m+2\overline{a})^2\, C_m(n)$$

$$+ (\overline{d}\,\overline{e} + \overline{b}\,\overline{d} + 2\,\overline{d}\,\overline{b}\,m + 2\,\overline{a}\,m^2\,\overline{b} + 2\,\overline{a}\,m\,\overline{b} + 2\,\overline{e}\,\overline{a}\,n - \overline{d}\,\overline{b}\,n)\,(\overline{d} + 2\,\overline{a}\,m + 2\,\overline{a})$$

$$(m+1)\,(\overline{d} + 3\,\overline{a} + 2\,\overline{a}\,m)\,(\overline{d} + \overline{a} + 2\,\overline{a}\,m)\,C_{m+1}(n) - (m+2)(-4\,\overline{a}^2\,\overline{c}\,m^2 \tag{33}$$

$$+ \overline{a}\,\overline{b}^2\,m^2 + 2\,\overline{a}\,\overline{b}^2\,m - 4\,\overline{a}\,\overline{c}\,m\,\overline{d} - 8\,\overline{a}^2\,\overline{c}\,m + m\,\overline{b}^2\,\overline{d} - \overline{a}\,\overline{e}^2 - \overline{d}^2\,\overline{c} + \overline{b}\,\overline{e}\,\overline{d} - 4\,\overline{a}^2\,\overline{c}$$

$$- 4\,\overline{a}\,\overline{c}\,\overline{d} + \overline{a}\,\overline{b}^2 + \overline{b}^2\,\overline{d})(\overline{a}\,m + \overline{a}\,n + \overline{a} + \overline{d})\,(m+1)\,(\overline{d} + 2\,\overline{a}\,m)\,C_{m+2}(n) = 0\,.$$

If $\overline{c} = 0$, then the recurrence equation

$$(n-m)\,(\overline{d}+2\,\overline{a}\,m)\,(\overline{d}+\overline{a}+2\,\overline{a}\,m)\,C_m(n) + (m+1)\,(\overline{b}\,m+\overline{e})\,(\overline{a}\,m+n\,\overline{a}+\overline{d})\,C_{m+1}(n) = 0 \tag{34}$$

is valid, and we get ($\overline{a}\,\overline{b} \neq 0$)

$$C_m(n) = \frac{\left(\frac{\overline{e}}{\overline{b}}\right)_n}{\left(\frac{\overline{d}}{\overline{a}}\right)_n} \left(-\frac{\overline{b}}{\overline{a}}\right)^n \cdot \frac{(-n)_m \left(\frac{\overline{d}}{2\overline{a}}\right)_m \left(\frac{\overline{a}+\overline{d}}{2\overline{a}}\right)_m}{\left(\frac{\overline{e}}{\overline{b}}\right)_m \left(\frac{\overline{a}n+\overline{d}}{\overline{a}}\right)_m m!} \left(\frac{4\overline{a}}{\overline{b}}\right)^m . \tag{35}$$

Therefore, the following representations for the powers in terms of the classical continuous orthogonal polynomials are valid:

$$(1-x)^n = 2^n\,\Gamma(\alpha+n+1) \sum_{m=0}^{n} \frac{(\alpha+\beta+2m+1)\,\Gamma(\alpha+\beta+m+1)}{\Gamma(\alpha+m+1)\,\Gamma(\alpha+\beta+n+m+2)}\,(-n)_m\,P_m^{(\alpha,\beta)}(x)$$

(see e.g. [20], 136, Eq. (2), or [17], § 5.2.4; note the essential misprint in this formula!),

$$(1+x)^n = 2^n\,\Gamma(\beta+n+1) \sum_{m=0}^{n} (-1)^m \frac{(\alpha+\beta+2m+1)\,\Gamma(\alpha+\beta+m+1)}{\Gamma(\beta+m+1)\,\Gamma(\alpha+\beta+n+m+2)}\,(-n)_m\,P_m^{(\alpha,\beta)}(x)\,,$$

$$x^n = \sum_{k=0}^{\lfloor n/2 \rfloor} \frac{(-n/2)_k\,(-n/2+1/2)_k\,(-n-\alpha)_k}{(-n/2-\alpha/2)_k\,(-n/2-\alpha/2+1/2)_k\,k!} \left(-\frac{1}{4}\right)^k \tilde{C}_{n-2k}^{\alpha}(x)\,,$$

$$\begin{aligned} x^n &= \frac{n!}{(\alpha)_n\,2^n} \sum_{k=0}^{\lfloor n/2 \rfloor} \frac{(-n/2-\alpha/2+1)_k\,(-n-\alpha)_k}{(-n/2-\alpha/2)_k\,k!}\,(-1)^k\,C_{n-2k}^{\alpha}(x) \\ &= \frac{n!}{2^n} \sum_{k=0}^{\lfloor n/2 \rfloor} \frac{n+\alpha-2k}{k!\,(\alpha)_{n+1-k}}\,C_{n-2k}^{\alpha}(x) \end{aligned}$$

(see e.g. [20], 144, Eq. (36), or [17], § 5.3.4),

$$x^n = (1+\alpha)_n \sum_{m=0}^{n} \frac{(-n)_m}{(1+\alpha)_m\,m!}\,(-1)^m\,\tilde{L}_m^{(\alpha)}(x)\,,$$

$$x^n = (1+\alpha)_n \sum_{m=0}^{n} \frac{(-n)_m}{(1+\alpha)_m}\,L_m^{(\alpha)}(x) = n! \sum_{m=0}^{n} \binom{n+\alpha}{n-m} (-1)^m\,L_m^{(\alpha)}(x)$$

(see e.g. [20], 118, Eq. (2), or [17], § 5.5.4),

$$x^n = \sum_{k=0}^{\lfloor n/2 \rfloor} \frac{(-n/2)_k\,(-n/2+1/2)_k}{k!}\,\tilde{H}_{n-2k}(x)\,,$$

$$x^n = \sum_{k=0}^{\lfloor n/2 \rfloor} \frac{(-n/2)_k\,(-n/2+1/2)_k}{k!\,2^{n-2k}} H_{n-2k}(x) = \frac{n!}{2^n} \sum_{k=0}^{\lfloor n/2 \rfloor} \frac{1}{k!\,(n-2k)!} H_{n-2k}(x)$$

(see e.g. [20], 110, Eq. (4), or [17], § 5.6.4),

$$x^n = \frac{(-2)^n}{(\alpha+2)_n} \sum_{m=0}^{n} \frac{(-n)_m\,(\alpha/2+1)_m\,(\alpha/2+3/2)_m}{(n+2+\alpha)_m\,m!}\, 2^m\, \widetilde{B}_m^{(\alpha)}(x)\ .$$

$$\begin{aligned} x^n &= \frac{(-2)^n}{(\alpha+2)_n} \sum_{m=0}^{n} \frac{(-n)_m\,(\alpha+1)_m\,(\alpha/2+3/2)_m}{(n+2+\alpha)_m\,(\alpha/2+1/2)_m\,m!} B_m^{(\alpha)}(x) \\ &= (-2)^n \sum_{m=0}^{n} (2m+\alpha+1)\, \frac{(-n)_m\,\Gamma(\alpha+m+1)}{m!\,\Gamma(n+m+\alpha+2)} B_m^{(\alpha)}(x) \end{aligned}$$

(see [2], (7.5); note the essential misprint in this formula!; compare [20], 150, Eq. (7)).

*Proof:* In § 2 it was shown how one obtains three essentially different cross rules for the connection coefficients between $P_n(x)$ and $Q_m(x)$. We modify this method here. For $Q_m(x)$, we have the differential equation

$$\overline{\sigma}(x)\, Q_m''(x) + \overline{\tau}(x)\, Q_m'(x) + \overline{\lambda}_m\, Q_m(x) = 0$$

with $\overline{\sigma}(x) = \overline{a}x^2 + \overline{b}x + \overline{c}$, and the derivative rule

$$\overline{\sigma}(x)\, Q_m'(x) = \overline{\alpha}_m\, Q_{m+1}(x) + \overline{\beta}_m\, Q_m(x) + \overline{\gamma}_m\, Q_{m-1}(x)\ ,$$

and it is easily seen that our current $P_n(x) = x^n$ satisfies any of the derivative rules

$$\overline{\sigma}(x)\, P_n'(x) = \overline{a}\, n\, P_{n+1}(x) + \overline{b}\, n\, P_n(x) + \overline{c}\, n\, P_{n-1}(x)\ . \tag{36}$$

Hence in our situation, we get the two cross rules (23) with $a_n = 1$, $b_n = c_n = 0$

$$C_m(n+1) = \overline{a}_{m-1}\, C_{m-1}(n) + \overline{b}_m\, C_m(n) + \overline{c}_{m+1}\, C_{m+1}(n) \tag{37}$$

and (25) with $\widehat{a}_n = 1/(n+1)$, $\widehat{b}_n = \widehat{c}_n = 0$

$$\frac{1}{n+1}\, C_m(n+1) = \overline{\widehat{a}}_{m-1}\, C_{m-1}(n) + \overline{\widehat{b}}_m\, C_m(n) + \overline{\widehat{c}}_{m+1}\, C_{m+1}(n) \tag{38}$$

which we had deduced in § 2. Using the derivative rule (36), we obtain the third cross rule

$$\overline{a}\, n\, C_m(n+1) + \overline{b}\, n\, C_m(n) + \overline{c}\, n\, C_m(n-1) = \overline{\alpha}_{m-1}\, C_{m-1}(n) + \overline{\beta}_m\, C_m(n) + \overline{\gamma}_{m+1}\, C_{m+1}(n)\ . \tag{39}$$

To receive the recurrence equation (33), we use Theorem 1 writing the cross rules in terms of $\overline{a}, \overline{b}, \overline{c}, \overline{d}$, and $\overline{e}$, only. Then by linear algebra we eliminate the variables $C_m(n+1)$ and $C_m(n-1)$ to obtain a pure recurrence equation w.r.t. $m$. (Similarly by elimination of the variables $C_{m-1}(n)$ and $C_{m+1}(n)$ a pure recurrence equation w.r.t. $n$ is obtained.) A shift by one gives (33).
If $\overline{c} = 0$, then the recurrence equation has still three terms, unfortunately. But since for $\overline{c} = 0$ in neither of the three cross rules (37)–(39) the variable $C_m(n-1)$ does occur, we can do a similar elimination, this time eliminating the variables $C_m(n+1)$ and $C_{m-1}(n)$, leading to the first order recurrence equation (34). Hence the hypergeometric representation (35) follows. The power representations for the Jacobi, Laguerre and Bessel polynomials are special cases thereof.
In the case of Hermite and Gegenbauer polynomials, (33) contains only the two terms $C_m(n)$ and $C_{m+2}(n)$, which leads to the desired representations. □

Note that, again, the recurrence equation (34) carries the complete information about the hypergeometric representations given in the theorem.
As an immediate consequence of the above theorem, we get the following connection coefficient results.

**Corollary 2** The following connection relations between the classical orthogonal polynomials are valid:

$$P_n^{(\alpha,\beta)}(x) = \sum_{m=0}^{n} (2m+\gamma+\beta+1)\frac{\Gamma(n+\beta+1)}{\Gamma(m+\beta+1)}\frac{\Gamma(n+m+\alpha+\beta+1)}{\Gamma(n+\alpha+\beta+1)} \cdot \frac{\Gamma(m+\gamma+\beta+1)}{\Gamma(n+m+\gamma+\beta+2)}\frac{(\alpha-\gamma)_{n-m}}{(n-m)!} P_m^{(\gamma,\beta)}(x) ,$$

(see e.g. [4], (13)),

$$P_n^{(\alpha,\beta)}(x) = \sum_{m=0}^{n} (-1)^{n-m}(2m+\alpha+\delta+1)\frac{\Gamma(n+\alpha+1)}{\Gamma(m+\alpha+1)}\frac{\Gamma(n+m+\alpha+\beta+1)}{\Gamma(n+\alpha+\beta+1)} \cdot \frac{\Gamma(m+\alpha+\delta+1)}{\Gamma(n+m+\alpha+\delta+2)}\frac{(\beta-\delta)_{n-m}}{(n-m)!} P_m^{(\alpha,\delta)}(x)$$

(see e.g. [6], (2.8)),

$$C_n^{\alpha}(x) = \frac{\Gamma(\beta)}{\Gamma(\alpha)\,\Gamma(\alpha-\beta)} \sum_{m=0}^{\lfloor n/2 \rfloor} \frac{(n-2k+\beta)\,\Gamma(k+\alpha-\beta)\,\Gamma(n-k+\alpha)}{k!\,\Gamma(n-k+\beta+1)} C_{n-2k}^{\beta}(x)$$

(see e.g. [5], (3.42)),

$$L_n^{(\alpha)}(x) = \sum_{m=0}^{n} \frac{(\alpha-\beta)_{n-m}}{(n-m)!} L_m^{(\beta)}(x)$$

(see e.g. [20], 119, Eq. (2)),

$$\begin{aligned} B_n^{(\alpha)}(x) &= \frac{(-1)^n\,(\alpha-\beta)_n}{(\beta+2)_n} \sum_{m=0}^{n} \frac{(-n)_m\,(\beta+1)_m\,(\beta/2+3/2)_m\,(n+\alpha+1)_m}{(n+2+\beta)_m\,(\beta/2+1/2)_m\,(\beta-\alpha+1-n)_m\,m!}(-1)^m B_m^{(\beta)}(x) \\ &= \sum_{m=0}^{n} (-1)^m\,(2m+\beta+1)\cdot\frac{(-n)_m\,\Gamma(\beta+m+1)\,(n+\alpha+1)_m\,\Gamma(\beta-\alpha+1)}{m!\,\Gamma(n+m+\beta+2)\,\Gamma(m-n+\beta-\alpha+1)} B_m^{(\beta)}(x) , \end{aligned}$$

(see [2], (8.2); note the essential misprint in this formula!).

*Proof:* We want to find the coefficients $C_m(n)$ in the relation (22)

$$P_n(x) = \sum_{m=0}^{n} C_m(n)\,Q_m(x) .$$

Combining

$$P_n(x) = \sum_{j\in\mathbb{Z}} A_j(n)\, x^j \qquad \text{and} \qquad x^j = \sum_{m\in\mathbb{Z}} B_m(j)\, Q_m(x)$$

yields the representation

$$P_n(x) = \sum_{j\in\mathbb{Z}} \sum_{m\in\mathbb{Z}} A_j(n)\, B_m(j)\, Q_m(x)\ ,$$

and interchanging the order of summation gives

$$C_m(n) = \sum_{j\in\mathbb{Z}} A_j(n)\, B_m(j)\ .$$

Similarly, if (as in the Gegenbauer case)

$$P_n(x) = \sum_{j\in\mathbb{Z}} A_j(n)\, x^{n-2j} \qquad \text{and} \qquad x^j = \sum_{m\in\mathbb{Z}} B_m(j)\, Q_{j-2m}(x)$$

one gets

$$D_m(n) = \sum_{j\in\mathbb{Z}} A_j(n)\, B_{m-j}(n-2j)$$

with

$$P_n(x) = \sum_{m=0}^{n} D_m(n)\, Q_{n-2m}(x)\ .$$

Since the summand $F(j,m,n) := A_j(n)\, B_m(j)$ turns out to be a hypergeometric term with respect to $(j,m,n)$, i.e., the term ratios $F(j+1,m,n)/F(j,m,n)$, $F(j,m+1,n)/F(j,m,n)$, and $F(j,m,n+1)/F(j,m,n)$ are rational functions, Zeilberger's algorithm ([26], [11], see e.g. [9]) applies and finds recurrence equations for $C_m(n)$ with respect to $m$ and $n$.
In all cases considered, Zeilberger's algorithm finds recurrence equations of first order with respect to $m$ (as well as for $n$). The given representations follow then from the initial value $C_n(n) = k_n/\overline{k}_n$. □

For some applications, it is important to know the rate of change in the direction of the parameters of the orthogonal systems, given in terms of the system itself. By a limiting process, these parameter derivative representations can be obtained from the results of Corollary 2.

**Corollary 3** The following representations for the parameter derivatives of the classical orthogonal polynomials are valid:

$$\begin{aligned}\frac{\partial}{\partial\alpha} P_n^{(\alpha,\beta)}(x) \quad &= \quad \sum_{m=0}^{n-1} \frac{1}{\alpha+\beta+1+m+n} \cdot \Big( P_n^{(\alpha,\beta)}(x) + \\ &\qquad \frac{\alpha+\beta+1+2\,m}{n-m}\, \frac{(\beta+m+1)_{n-m}}{(\alpha+\beta+m+1)_{n-m}}\, P_m^{(\alpha,\beta)}(x) \Big)\end{aligned}$$

(see [7], Theorem 3),

$$\frac{\partial}{\partial\alpha} \widetilde{P}_n^{(\alpha,\beta)}(x) = \sum_{m=0}^{n-1} \frac{2^{n-m}}{n-m} \frac{\binom{2m+\alpha+\beta}{m}}{\binom{2n+\alpha+\beta}{n}} \frac{\alpha+\beta+1+2\,m}{\alpha+\beta+1+m+n} \frac{(\beta+m+1)_{n-m}}{(\alpha+\beta+m+1)_{n-m}} \widetilde{P}_m^{(\alpha,\beta)}(x)\ ,$$

$$\frac{\partial}{\partial\beta} P_n^{(\alpha,\beta)}(x) = \sum_{m=0}^{n-1} \frac{1}{\alpha+\beta+1+m+n} \cdot \Big( P_n^{(\alpha,\beta)}(x) + (-1)^{n-m} \frac{\alpha+\beta+1+2m}{n-m} \frac{(\alpha+m+1)_{n-m}}{(\alpha+\beta+m+1)_{n-m}} P_m^{(\alpha,\beta)}(x) \Big)$$

(see [7], Theorem 3),

$$\frac{\partial}{\partial\beta} \widetilde{P}_n^{(\alpha,\beta)}(x) = \sum_{m=0}^{n-1} \frac{(-2)^{n-m}}{n-m} \frac{\binom{2m+\alpha+\beta}{m}}{\binom{2n+\alpha+\beta}{n}} \frac{\alpha+\beta+1+2m}{\alpha+\beta+1+m+n} \frac{(\alpha+m+1)_{n-m}}{(\alpha+\beta+m+1)_{n-m}} \widetilde{P}_m^{(\alpha,\beta)}(x) \ ,$$

$$\begin{aligned} \frac{\partial}{\partial\alpha} C_n^{\alpha}(x) &= \sum_{m=0}^{n-1} \left( \frac{2\,(1+m)}{(2\alpha+m)\,(2\alpha+1+2m)} + \frac{2}{2\,\alpha+m+n} \right) C_n^{\alpha}(x) \\ &\quad + \sum_{m=0}^{n-1} \frac{2\,(1+(-1)^{n-m})\,(\alpha+m)}{(2\alpha+m+n)\,(n-m)} C_m^{\alpha}(x) \end{aligned}$$

(see [12], Theorem 10),

$$\begin{aligned} \frac{\partial}{\partial\alpha} \widetilde{C}_n^{\alpha}(x) &= \sum_{m=0}^{n-1} 2^{m-n+1} \frac{(\alpha)_m\, n!}{(\alpha)_n\, m!} \frac{(1+(-1)^{n-m})\,(\alpha+m)}{(2\alpha+m+n)\,(n-m)} \widetilde{C}_m^{\alpha}(x) \\ &= \sum_{k=1}^{\lfloor n/2 \rfloor} \frac{n!}{(\alpha+n-2k)_{2k}\, 4^k\, (n-2k)!} \frac{n-2k+\alpha}{k\,(n-k+\alpha)} \widetilde{C}_{n-2k}^{\alpha}(x) \ , \end{aligned}$$

$$\frac{\partial}{\partial\alpha} L_n^{(\alpha)}(x) = \sum_{m=0}^{n-1} \frac{1}{n-m} L_m^{(\alpha)}(x) \ ,$$

(see [12], Theorem 10),

$$\frac{\partial}{\partial\alpha} \widetilde{L}_n^{(\alpha)}(x) = \sum_{m=0}^{n-1} \frac{(-1)^{n-m}}{n-m} \frac{n!}{m!} \widetilde{L}_m^{(\alpha)}(x) \ ,$$

$$\frac{\partial}{\partial\alpha} B_n^{\alpha}(x) = \sum_{m=0}^{n-1} \frac{1}{\alpha+n+m+1} \cdot \Big( B_n^{\alpha}(x) + (-1)^{n-m} \frac{2m+\alpha+1}{(n-m)} \frac{n!}{(\alpha+m+1)_{n-m}\, m!} B_m^{\alpha}(x) \Big) \ ,$$

$$\frac{\partial}{\partial\alpha} \widetilde{B}_n^{\alpha}(x) = \sum_{m=0}^{n-1} \frac{(-2)^{n-m}}{(n-m)\,(\alpha+n+m+1)} \frac{n!}{(\alpha+2m+2)_{2n-2m-1}\, m!} \widetilde{B}_m^{\alpha}(x) \ .$$

*Proof:* Given the connection relation

$$P_n^{\alpha}(x) = \sum_{m=0}^{n} C_m(n;\alpha,\beta)\, P_m^{\beta}(x) \ ,$$

we build the difference quotient

$$\begin{aligned}\frac{P_n^\alpha(x)-P_n^\beta(x)}{\alpha-\beta} &= \sum_{m=0}^{n}\frac{C_m(n;\alpha,\beta)}{\alpha-\beta}\,P_m^\beta(x)-\frac{P_n^\beta(x)}{\alpha-\beta}\\ &= \frac{C_n(n;\alpha,\beta)-1}{\alpha-\beta}\,P_n^\beta(x)+\sum_{m=0}^{n-1}\frac{C_m(n;\alpha,\beta)}{\alpha-\beta}\,P_m^\beta(x)\ ,\end{aligned}$$

so that with $\beta\to\alpha$

$$\frac{\partial}{\partial\alpha}P_n^\alpha(x)=\lim_{\beta\to\alpha}\frac{C_n(n;\alpha,\beta)-1}{\alpha-\beta}\,P_n^\alpha(x)+\sum_{m=0}^{n-1}\lim_{\beta\to\alpha}\frac{C_m(n;\alpha,\beta)}{\alpha-\beta}\,P_m^\alpha(x) \tag{40}$$

since the systems $P_n^\alpha(x)$ are continuous w.r.t. $\alpha$. This gives the results. □

Note that for monic polynomials (and moreover if $k_n$ does not depend on $\alpha$ as in the Laguerre case) the first limit in (40) equals zero. Hence the parameter derivative representations are simplest in such a case.

## 5 Hypergeometric Representations: Discrete Case

By $h_n^{(\alpha,\beta)}(x,N)$ and $Q_n(x;\alpha,\beta,N)$ we denote two commonly used standardizations of the *Hahn polynomials* (see [19], and [23]), and by $m_n^{(\gamma,\mu)}(x)$, $k_n^{(p)}(x,N)$ and $c_n^{(\mu)}(x)$ the *Meixner, Krawchouk* and *Charlier polynomials* are denoted, respectively. They have the standardizations

| system | $h_n^{(\alpha,\beta)}(x,N)$ | $Q_n(x;\alpha,\beta,N)$ | $m_n^{(\gamma,\mu)}(x)$ | $k_n^{(p)}(x,N)$ | $c_n^{(\mu)}(x)$ | $K_n^{(\alpha,\beta)}(x)$ |
|---|---|---|---|---|---|---|
| $k_n$ | $\binom{\alpha+\beta+2n}{n}$ | $\frac{(\alpha+\beta+n+1)_n}{(-N)_n\,(\alpha+1)_n}$ | $\left(1-\frac{1}{\mu}\right)^n$ | $\frac{1}{n!}$ | $\left(-\frac{1}{\mu}\right)^n$ | $\alpha^n$ |

The polynomials $t_n(x,N):=h_n^{(0,0)}(x,N)$ are the *discrete Chebyshev polynomials.* The polynomials $K_n^{(\alpha,\beta)}(x)$ given in (21), are not orthogonal, but satisfy the difference equation

$$\Delta\nabla y(x)+(\alpha x+\beta)\,\Delta y(x)+\lambda_n y(x)=0\ .$$

The monic counterparts of the discrete systems will be denoted by $\widetilde{h}_n^{(\alpha,\beta)}(x,N)$, $\widetilde{Q}_n(x;\alpha,\beta,N)$, $\widetilde{t}_n(x,N)$, $\widetilde{m}_n^{(\gamma,\mu)}(x)$, $\widetilde{k}_n^{(p)}(x,N)$ and $\widetilde{c}_n^{(\mu)}(x)$, respectively. Observe that therefore by $\widetilde{h}_n$ we do *not* denote the Hahn-Eberlein polynomials $\widetilde{h}_n^{(\mu,\nu)}(x,N)$ as in [19].
In the continuous case, we looked for power series representations, i.e., we set $Q_m(x)=x^m$. The corresponding choice in the discrete case is a representation in terms of the falling factorial

$$Q_m(x)=x^{\underline{m}}:=x\,(x-1)\cdots(x-m+1)=(x-m+1)_m=(-1)^m\,(-x)_m\ .$$

We get

**Theorem 6** Let $P_n(x)$ be a monic polynomial system given by the difference equation (4) with $\sigma(x)=ax^2+bx+c$, and $\tau(x)=dx+e$. Then the series coefficients $C_m(n)$ given by

$$P_n(x)=\sum_{m=0}^{n}C_m(n)\,x^{\underline{m}} \tag{41}$$

satisfy the recurrence equation

$$\begin{aligned}&(\,a\,n+a\,m-a+d\,)\,(\,n-m\,)\,C_m(n)\\&+(\,m+1\,)\,(\,a\,n^2-2\,a\,m^2-a\,n-a\,m+n\,d-2\,d\,m-b\,m-d-e\,)\,C_{m+1}(n)\\&-(\,m+1\,)\,(\,m+2\,)\,(\,a\,m^2+2\,a\,m+d\,m+b\,m+a+d+b+c+e\,)\,C_{m+2}(n)=0\;.\end{aligned} \tag{42}$$

If $c=0$, then the recurrence equation

$$(\,n-m\,)\,(\,a\,m+d+a\,n-a\,)\,C_m(n)-(\,m+1\,)\,(\,a\,m^2+m\,b+m\,d+e\,)\,C_{m+1}(n)=0 \tag{43}$$

is valid, and we have the hypergeometric representation

$$\begin{aligned}P_n(x)\;&=\;\frac{\left(\frac{d}{a}-1\right)_n\left(\frac{b+d+\sqrt{(b+d)^2-4ae}}{2a}\right)_n\left(\frac{b+d-\sqrt{(b+d)^2-4ae}}{2a}\right)_n}{\left(\frac{d}{2a}\right)_n\left(\frac{d-a}{2a}\right)_n}\left(\frac{1}{4}\right)^n\\&\cdot\,{}_3F_2\left(\left.\begin{matrix}-n,-x,n-1+\frac{d}{a}\\ \frac{b+d+\sqrt{(b-d)^2-4ae}}{2a},\frac{b+d-\sqrt{(b-d)^2-4ae}}{2a}\end{matrix}\right|1\right)\;,\end{aligned} \tag{44}$$

valid for $a\neq 0$, or

$$P_n(x)=\left(\frac{e}{b+d}\right)_n\left(1+\frac{b}{d}\right)^n\,{}_2F_1\left(\left.\begin{matrix}-n,-x\\ \frac{e}{b+d}\end{matrix}\right|\frac{d}{b+d}\right)\;, \tag{45}$$

valid for $a=0, b+d\neq 0$, or finally

$$P_n(x)=\left(-\frac{e}{b}\right)^n\,{}_2F_0\left(\left.\begin{matrix}-n,-x\\ -\end{matrix}\right|-\frac{b}{e}\right)\;, \tag{46}$$

valid for $a=0, d=-b$.
Therefore, the classical discrete orthogonal polynomials and their monic counterparts have the following hypergeometric series representations:

$$h_n^{(\alpha,\beta)}(x,N)=\frac{(-1)^n}{n!}\,(\beta+1)_n\,(N-n)_n\,{}_3F_2\left(\left.\begin{matrix}-n,-x,n+1+\alpha+\beta\\ \beta+1,1-N\end{matrix}\right|1\right)$$

(see e.g. [19], p. 54, Table 2.4),

$$\tilde{h}_n^{(\alpha,\beta)}(x,N)=\frac{(1+\beta)_n\,(1-N)_n}{(1+n+\alpha+\beta)_n}\,{}_3F_2\left(\left.\begin{matrix}-n,-x,n+1+\alpha+\beta\\ \beta+1,1-N\end{matrix}\right|1\right)\;,$$

$$t_n(x,N)=(-1)^n\,(N-n)_n\cdot{}_3F_2\left(\left.\begin{matrix}-n,-x,n+1\\ 1,1-N\end{matrix}\right|1\right)\;,$$

$$\tilde{t}_n(x,N)=\frac{n!\,(1-N)_n}{(1/2)_n\,4^n}\,{}_3F_2\left(\left.\begin{matrix}-n,-x,n+1\\ 1,1-N\end{matrix}\right|1\right)\;,$$

$$Q_n(x;\alpha,\beta,N)={}_3F_2\left(\left.\begin{matrix}-n,-x,n+1+\alpha+\beta\\ \alpha+1,-N\end{matrix}\right|1\right)$$

(see e.g. [23], 1.5)

$$\tilde{Q}_n(x;\alpha,\beta,N) = \tilde{h}_n^{(\beta,\alpha)}(x,N+1) \;,$$

$$m_n^{(\gamma,\mu)}(x) = (\gamma)_n \; {}_2F_1\left(\begin{array}{c} -n,-x \\ \gamma \end{array} \middle| 1-\frac{1}{\mu}\right)$$

(see e.g. [19], p. 54, Table 2.4),

$$\tilde{m}_n^{(\gamma,\mu)}(x) = (\gamma)_n \; \left(\frac{\mu}{\mu-1}\right)^n \; {}_2F_1\left(\begin{array}{c} -n,-x \\ \gamma \end{array} \middle| 1-\frac{1}{\mu}\right) \;,$$

$$k_n^{(p)}(x,N) = (-1)^n \binom{N}{n} p^n \; {}_2F_1\left(\begin{array}{c} -n,-x \\ -N \end{array} \middle| \frac{1}{p}\right)$$

(see e.g. [19], p. 54, Table 2.4),

$$\tilde{k}_n^{(p)}(x,N) = (-N)_n \, p^n \; {}_2F_1\left(\begin{array}{c} -n,-x \\ -N \end{array} \middle| \frac{1}{p}\right) \;,$$

$$c_n^{(\mu)}(x) = {}_2F_0\left(\begin{array}{c} -n,-x \\ - \end{array} \middle| -\frac{1}{\mu}\right)$$

(see e.g. [19], p. 54, Table 2.4),

$$\tilde{c}_n^{(\mu)}(x) = (-\mu)^n \; {}_2F_0\left(\begin{array}{c} -n,-x \\ - \end{array} \middle| -\frac{1}{\mu}\right) \,.$$

These results are all particular cases of the recurrence equation (43).

*Proof:* Substituting the series (41) into the difference equation, and equating the coefficients of the falling factorials yields the recurrence equation (42) which had been obtained by Lesky [16].

This conversion can be easily done using a computer algebra system by bringing the given difference equation into the form (5), expanding it, and replacing any occurrence of $\Delta y(x)$ by $(m+1)\,C_{m+1}$, any occurrence of a product $x\,y(x)$ by $C_{m-1} + m\,C_m$ and any occurrence of a shift $y(x+1)$ by $C_m + (m+1)\,C_{m+1}$ since

$$\Delta x^{\underline{m}} = m\,x^{\underline{m-1}} \;, \quad x\,x^{\underline{m}} = x^{\underline{m+1}} + m\,x^{\underline{m}} \;, \quad \text{and} \quad (x+1)^{\underline{m}} = x^{\underline{m}} + m\,x^{\underline{m-1}} \,.$$

Iteratively for all nonnegative integers $j,k$ any of the terms $x^j\,\Delta^k y(x)$ and $x^j\,y(x+k)$ can be replaced by these rules. Note that this method can also be applied for higher order difference equations with polynomials coefficients.

Different from the continuous case, the recurrence equation (42) does *not* degenerate to a two-term recurrence equation for $c = 0$. To get (43), nevertheless, we must use a different approach. One possibility is to apply Petkovšek's algorithm to the recurrence equation (42), leading to (43).

Another possibility is to modify the method which will be used in the next section to deduce representations of the falling factorials in terms of discrete orthogonal systems. This method yields (43) directly.

As soon as (43) is deduced, the initial value $C_n(n) = 1$ gives the hypergeometric representations (44)–(46) which include all other representations by substituting the particular parameters $a, b, c, d$, and $e$. $\square$

We would like to mention that, again, a single recurrence equation, (42), carries the complete information about the hypergeometric representations given in the theorem.
Note furthermore, that the radicals in (44) do only occur by the representation used: the radical factors come in pairs whose product is radical-free. Note that the computation which gives (44), answers a question raised by Koornwinder [15]. For more examples of this type see [13].
Our method was able to find hypergeometric series representations for the particular case $c = 0$. This is the most important situation since all the classical discrete orthogonal families are of this type, corresponding to the fact that their discrete support has zero as left boundary point (see e.g. [19], Tables 2.1–2.3).
By construction, all the series representations determined have an upper parameter $-x$. The question remains, however, whether or not such a hypergeometric series representation might be valid for $c \neq 0$, too. In general, the answer is no. Petkovšek's algorithm shows that the recurrence equation (42) does *not* generally have a hypergeometric term solution.
Note that the hypergeometric representation (21) for $K_n^{(\alpha,\beta)}(x)$ is not of this type, and cannot be obtained by the given method. By Petkovšek's algorithm there is no representation (41) with a hypergeometric term $C_m(n)$ for these polynomials.

## 6 Falling Factorial Representations

Whereas in the last section we considered the specific connection coefficient problem for $Q_m(x) = x^{\underline{m}}$, in this section the opposite problem, having $P_n(x) = x^{\underline{n}}$, is studied.

**Theorem 7** Let $Q_m(x)$ be a monic polynomial system given by the difference equation (4) with $\overline{\sigma}(x) = \overline{a}x^2 + \overline{b}x + \overline{c}$, and $\overline{\tau}(x) = \overline{d}x + \overline{e}$. Then the coefficients $C_m(n)$ of the falling factorial representations

$$x^{\underline{n}} = \sum_{m=0}^{n} C_m(n)\, Q_m(x) \tag{47}$$

satisfy the recurrence equation

$$\begin{aligned}
&(2m\overline{a}+\overline{a}+\overline{d})(2m\overline{a}+3\overline{a}+\overline{d})(2m\overline{a}+2\overline{a}+\overline{d})^2(2m\overline{a}+\overline{d})(n-m)\,C_m(n)+\\
&(2m\overline{a}+\overline{a}+\overline{d})(2m\overline{a}+3\overline{a}+\overline{d})(2m\overline{a}+2\overline{a}+\overline{d})(m+1)(2m^2n\overline{a}^2-2m^2\overline{a}^2\\
&+m^2\overline{a}\,\overline{d}+2m^2\overline{a}\,\overline{b}+2mn\overline{a}^2+2mn\overline{a}\,\overline{d}-2m\overline{a}^2-m\overline{a}\,\overline{d}+2m\overline{a}\,\overline{b}+m\overline{d}^2\\
&+2m\overline{d}\,\overline{b}+n\overline{a}\,\overline{d}+2n\overline{a}\,\overline{e}-n\overline{d}\,\overline{b}-\overline{a}\,\overline{d}+\overline{d}\,\overline{b}+\overline{d}\,\overline{e})\,C_{m+1}(n)\\
&+(m+1)(2m\overline{a}+\overline{d})(m^4\overline{a}^3+4m^3\overline{a}^3+2m^3\overline{a}^2\overline{d}+6m^2\overline{a}^3+6m^2\overline{a}^2\overline{d}\\
&+4m^2\overline{a}^2\overline{c}+2m^2\overline{a}^2\overline{e}+m^2\overline{a}\,\overline{d}^2-m^2\overline{a}\,\overline{d}\,\overline{b}-m^2\overline{a}\,\overline{b}^2+4m\overline{a}^3+6m\overline{a}^2\overline{d}\\
&+8m\overline{a}^2\overline{c}+4m\overline{a}^2\overline{e}+2m\overline{a}\,\overline{d}^2-2m\overline{a}\,\overline{d}\,\overline{b}+4m\overline{a}\,\overline{d}\,\overline{c}+2m\overline{a}\,\overline{d}\,\overline{e}\\
&-2m\overline{a}\,\overline{b}^2-m\overline{d}^2\overline{b}-m\overline{d}\,\overline{b}^2+\overline{a}^3+2\overline{a}^2\overline{d}+4\overline{a}^2\overline{c}+2\overline{a}^2\overline{e}+\overline{a}\,\overline{d}^2\\
&-\overline{a}\,\overline{d}\,\overline{b}+4\overline{a}\,\overline{d}\,\overline{c}+2\overline{a}\,\overline{d}\,\overline{e}-\overline{a}\,\overline{b}^2+\overline{a}\,\overline{e}^2-\overline{d}^2\overline{b}+\overline{d}^2\overline{c}\\
&-\overline{d}\,\overline{b}^2-\overline{d}\,\overline{b}\,\overline{e})(m+2)(m\overline{a}+n\overline{a}+\overline{a}+\overline{d})\,C_{m+2}(n)=0\,.
\end{aligned} \tag{48}$$

If $\overline{c}=0$, then

$$0 = (\overline{d}+\overline{a}+2\,\overline{a}\,m)\,(\overline{d}+2\,\overline{a}\,m)\,(-n+m)\,C_m(n) \\ -(\overline{a}\,n+\overline{d}+\overline{a}\,m)\,(m+1)\,(\overline{a}\,m^2+m\,\overline{d}+m\,\overline{b}+\overline{e})\,C_{m+1}(n)\ . \tag{49}$$

Therefore, the following representations for the falling factorials in terms of the classical discrete orthogonal polynomials are valid:

$$x^{\underline{n}} = \sum_{m=0}^{n} \frac{(\beta+1)_n\,(1-N)_n\,(-1)^n\,(1+\alpha+\beta+2m)\,(-n)_m\,(1+\alpha+\beta)_m}{(\alpha+\beta+2)_n\,(1+\alpha+\beta)\,(n+2+\alpha+\beta)_m\,(\beta+1)_m\,(1-N)_m}\,h_m^{(\alpha,\beta)}(x,N)\ ,$$

$$x^{\underline{n}} = \sum_{m=0}^{n} \frac{(\beta+1)_n\,(1-N)_n\,(-1)^n}{(\alpha+\beta+2)_n}\,\frac{(-n)_m\,(\alpha/2+\beta/2+1)_m\,(\alpha/2+\beta/2+3/2)_m\,4^m}{(n+2+\alpha+\beta)_m\,(\beta+1)_m\,(1-N)_m\,m!}\,\tilde{h}_m^{(\alpha,\beta)}(x,N)\ ,$$

$$x^{\underline{n}} = \sum_{m=0}^{n} \frac{(1+\alpha)_n\,(-N)_n\,(-1)^n}{(\alpha+\beta+2)_n}\,\frac{(\alpha+\beta+1+2m)}{(\alpha+\beta+1)}\,\frac{(-n)_m\,(1+\alpha+\beta)_m}{(n+2+\alpha+\beta)_m\,m!}\,Q_m(x;\alpha,\beta,N)$$

(compare [8], (4.2)–(4.3)),

$$x^{\underline{n}} = \frac{(1-N)_n\,(-1)^n}{n+1}\,\sum_{m=0}^{n} \frac{(-n)_m\,(1+2m)}{(n+2)_m\,(1-N)_m}\,t_m(x,N)\ ,$$

$$x^{\underline{n}} = \frac{(1-N)_n\,(-1)^n}{n+1}\,\sum_{m=0}^{n} \frac{(-n)_m\,(3/2)_m\,4^m}{(n+2)_m\,(1-N)_m\,m!}\,\tilde{t}_m(x,N)\ ,$$

$$x^{\underline{n}} = \sum_{m=0}^{n} \frac{(-1)^n\,(\gamma)_n\,\left(\frac{\mu}{\mu-1}\right)^n\,(-n)_m}{(\gamma)_m\,m!}\,m_m^{(\gamma,\mu)}(x)\ ,$$

$$x^{\underline{n}} = \sum_{m=0}^{n} \frac{(-1)^n\,(\gamma)_n\,\left(\frac{\mu}{\mu-1}\right)^{n-m}\,(-n)_m}{(\gamma)_m\,m!}\,\tilde{m}_m^{(\gamma,\mu)}(x)\ ,$$

$$x^{\underline{n}} = \sum_{m=0}^{n} \frac{(-1)^n\,(-N)_n\,p^{n-m}\,(-n)_m}{(-N)_m}\,k_m^{(p)}(x,N)\ ,$$

$$x^{\underline{n}} = \sum_{m=0}^{n} \frac{(-1)^n\,(-N)_n\,p^{n-m}\,(-n)_m}{(-N)_m\,m!}\,\tilde{k}_m^{(p)}(x,N)\ ,$$

$$x^{\underline{n}} = \sum_{m=0}^{n} \frac{\mu^n\,(-n)_m}{m!}\,c_m^{(\mu)}(x)\ ,$$

$$x^{\underline{n}} = \sum_{m=0}^{n} \frac{(-1)^n\,(-\mu)^{n-m}\,(-n)_m}{m!}\,\tilde{c}_m^{(\mu)}(x)\ .$$

*Proof:* In § 2 it was shown how one obtains three essentially different cross rules for the connection coefficients between $P_n(x)$ and $Q_m(x)$. We modify this method here. For $Q_m(x)$, we have the difference equation

$$\overline{\sigma}(x)\,\Delta\nabla Q_m(x)+\overline{\tau}(x)\,\Delta Q_m(x)+\overline{\lambda}_m\,Q_m(x)=0$$

with $\overline{\sigma}(x) = \overline{a}x^2 + \overline{b}x + \overline{c}$, and the difference rule (11)

$$(\overline{\sigma}(x) + \overline{\tau}(x))\,\Delta Q_m(x) = \overline{\alpha}_m\,Q_{m+1}(x) + (\overline{\beta}_m - \overline{\lambda}_m)\,Q_m(x) + \overline{\gamma}_m\,Q_{m-1}(x)\ ,$$

and it is easily seen that our current $P_n(x) = x^{\underline{n}}$ satisfies any of the difference rules

$$\begin{aligned}(\overline{\sigma}(x) + \overline{\tau}(x))\,\Delta P_n(x) \quad &= \quad \overline{a}\,n\,P_{n+1}(x) + n\left(\overline{a}\,(2n-1) + \overline{b} + \overline{d}\right) P_n(x) \\ &\quad + n\left((n-1)(\overline{a}(n-1) + \overline{b} + \overline{d}) + \overline{c} + \overline{e}\right) P_{n-1}(x)\ .\end{aligned} \tag{50}$$

Hence in our situation, we get the two cross rules (23) with $a_n = 1$, $b_n = n$, $c_n = 0$

$$C_m(n+1) = \overline{a}_{m-1}\,C_{m-1}(n) + \overline{b}_m\,C_m(n) + \overline{c}_{m+1}\,C_{m+1}(n) \tag{51}$$

and (25) with $\widehat{a}_n = 1/(n+1)$, $\widehat{b}_n = \widehat{c}_n = 0$

$$\frac{1}{n+1}a_n\,C_m(n+1) = \overline{\overline{a}}_{m-1}\,C_{m-1}(n) + \overline{\overline{b}}_m\,C_m(n) + \overline{\overline{c}}_{m+1}\,C_{m+1}(n) \tag{52}$$

which we had deduced in § 2. Using the difference rule (50), we obtain the third cross rule

$$\begin{aligned}\overline{a}\,n\,C_m(n+1) + n\left(\overline{a}\,(2n-1) + \overline{b} + \overline{d}\right) C_m(n) + n\left((n-1)(\overline{a}(n-1) + \overline{b} + \overline{d}) + \overline{c} + \overline{e}\right) C_m(n-1) \\ = \overline{\alpha}_{m-1}\,C_{m-1}(n) + \overline{\beta}_m\,C_m(n) + \overline{\gamma}_{m+1}\,C_{m+1}(n)\,.\end{aligned} \tag{53}$$

To receive the recurrence equation (48), we use Theorem 1 writing the cross rules in terms of $\overline{a}, \overline{b}, \overline{c}, \overline{d}$, and $\overline{e}$, only. Then by linear algebra we eliminate the variables $C_m(n+1)$ and $C_m(n-1)$ to obtain a pure recurrence equation w.r.t. $m$. (Similarly by elimination of the variables $C_{m-1}(n)$ and $C_{m+1}(n)$ a pure recurrence equation w.r.t. $n$ is obtained.) A shift by one gives (48).
If $\overline{c} = 0$, then the recurrence equation has still three terms, unfortunately. For $\overline{c} = 0$, we find a fourth cross rule to eliminate one more variable in the following way. Since the (second) difference rule (10)

$$\overline{\sigma}(x)\,\nabla Q_m(x) = \overline{\alpha}_m\,Q_{m+1}(x) + \overline{\beta}_m\,Q_m(x) + \overline{\gamma}_m\,Q_{m-1}(x)$$

is valid, we can use the fact that for $\overline{c} = 0$ any of the difference rules

$$\overline{\sigma}(x)\,\nabla P_n(x) = \overline{a}\,n\,P_{n+1}(x) + n\,(\overline{a}\,n + \overline{b})\,P_n(x)\ ,$$

and therefore the fourth cross rule

$$\overline{a}\,n\,C_m(n+1) + n\,(\overline{a}\,n + \overline{b})\,C_m(n) = \overline{\alpha}_{m-1}\,C_{m-1}(n) + \overline{\beta}_m\,C_m(n) + \overline{\gamma}_{m+1}\,C_{m+1}(n) \tag{54}$$

is valid.
Eliminating the variables $C_m(n+1)$, $C_{m-1}(n)$, and $C_{m-1}(n)$ from the four cross rules (51)–(54) gives the first order recurrence equation (49). This leads to the desired representations. □

Whereas we admit that all the shifted factorial representations of the theorem essentially were known ([8], [22]), our presentation unifies this development. In particular, the recurrence equation (49) carries the complete information about the falling factorial representations given in the theorem.
Petkovšek's algorithm proves, again, that for the family $K_n^{(\alpha,\beta)}(x)$ there is no representation (47) with a hypergeometric term $C_m(n)$.
As an immediate consequence of Theorem 7, we get the following connection coefficient results.

**Corollary 4** The following connection relations between the classical discrete orthogonal polynomials are valid:

$$h_n^{(\alpha,\beta)}(x,N) = \sum_{m=0}^{n} \frac{(\beta-\delta)_n\,(1-N)_n\,(\alpha+1)_n}{(2+\alpha+\delta)_n\,n!}\cdot$$

$$\cdot\,\frac{(\alpha+\delta+1+2m)}{(\alpha+\delta+1)}\,\frac{(-n)_m\,(1+\alpha+\delta)_m\,(n+1+\alpha+\beta)_m}{(1-N)_m\,(\alpha+1)_m\,(\alpha+2+n+\delta)_m\,(-n-\beta+\delta+1)_m}\,h_m^{(\alpha,\delta)}(x,N)\;, \tag{55}$$

$$\tilde{h}_n^{(\alpha,\beta)}(x,N) = \sum_{m=0}^{n} \frac{(\alpha+1)_n\,(1-N)_n\,(1+\alpha+\beta)_n\,(\beta-\delta)_n}{(2+\alpha+\delta)_n\,(\alpha/2+\beta/2+1/2)_n\,(\alpha/2+\beta/2+1)_n\,4^n}\cdot$$

$$\cdot\,\frac{(-n)_m\,(n+1+\alpha+\beta)_m\,(\alpha/2+\delta/2+1)_m\,(\alpha/2+\delta/2+3/2)_m\,4^m}{(1-N)_m\,(\alpha+1)_m\,(\alpha+2+n+\delta)_m\,(-n-\beta+\delta+1)_m\,m!}\,\tilde{h}_m^{(\alpha,\delta)}(x,N)\;,$$

$$h_n^{(\alpha,\beta)}(x,N) = \sum_{m=0}^{n} \frac{(\beta+1)_n\,(1-N)_n\,(\alpha-\gamma)_n\,(-1)^n}{(2+\beta+\gamma)_n\,n!}\cdot$$

$$\cdot\,\frac{(\beta+\gamma+1+2m)}{(\beta+\gamma+1)}\,\frac{(-n)_m\,(1+\beta+\gamma)_m\,(n+1+\alpha+\beta)_m\,(-1)^m}{(1-N)_m\,(\beta+1)_m\,(\beta+\gamma+n+2)_m\,(\gamma-\alpha-n+1)_m}h_m^{(\gamma,\beta)}(x,N)\;,$$

$$\tilde{h}_n^{(\alpha,\beta)}(x,N) = \sum_{m=0}^{n} \frac{(\beta+1)_n\,(1-N)_n\,(1+\alpha+\beta)_n\,(\alpha-\gamma)_n\,(-1)^n}{(2+\beta+\gamma)_n\,(\alpha/2+\beta/2+1/2)_n\,(\alpha/2+\beta/2+1)_n\,4^n}\cdot$$

$$\cdot\,\frac{(-n)_m\,(n+1+\alpha+\beta)_m\,(\beta/2+\gamma/2+1)_m\,(\beta/2+\gamma/2+3/2)_m\,(-4)^m}{(\beta+\gamma+n+2)_m\,(1-N)_m\,(\beta+1)_m\,(\gamma-\alpha-n+1)_m\,m!}\,\tilde{h}_m^{(\gamma,\beta)}(x,N)\;,$$

$$\tilde{h}_n^{(\alpha,\alpha)}(x,N) = \tag{56}$$

$$\sum_{k=0}^{\lfloor n/2\rfloor} \frac{(-n/2)_k\,(-(n-1)/2)_k\,(\alpha-\gamma)_k\,((N-n)/2)_k\,(N-n+1)/2)_k\,(-n-\gamma-1/2)_k}{(1/4-\gamma/2-n/2)_k\,(-n+1/2-\alpha)_k\,(-n/2-1/4-\gamma/2)_k\,k!\,4^k}\,\tilde{h}_{n-2k}^{(\gamma,\gamma)}(x,N)\;,$$

$$h_n^{(\alpha,\alpha)}(x,N) = \frac{(\alpha+1)_n\,(\alpha+1/2)_n\,(2\gamma+1)_n}{(\gamma+1)_n\,(\gamma+1/2)_n\,(2\alpha+1)_n}\cdot \tag{57}$$

$$\sum_{k=0}^{\lfloor n/2\rfloor} \frac{\left(\frac{N-n}{2}\right)_k\left(\frac{N-n+1}{2}\right)_k(\alpha-\gamma)_k\left(\frac{3-2\gamma-2n}{4}\right)_k\left(\frac{-2\gamma-2n-1}{2}\right)_k\left(\frac{-\gamma-n}{2}\right)_k\left(\frac{-\gamma-n+1}{2}\right)_k 4^k}{(-\gamma-n/2)_k\,(-\gamma/2-n/2-1/4)_k\,(-n-\alpha+1/2)_k\,(-\gamma-n/2+1/2)_k\,k!}\,h_{n-2k}^{(\gamma,\gamma)}(x,N)\;,$$

$$Q_n(x;\alpha,\alpha,N) = \frac{(\alpha+1/2)_n\,(2\gamma+1)_n}{(\gamma+1/2)_n\,(2\alpha+1)_n}\cdot \tag{58}$$

$$\sum_{k=0}^{\lfloor n/2\rfloor} \frac{(-n/2)_k\,(-(n-1)/2)_k\,(\alpha-\gamma)_k\,(3/4-\gamma/2-n/2)_k\,(-\gamma-n-1/2)_k}{(-\gamma-n/2)_k\,(-\gamma/2-n/2-1/4)_k\,(-n-\alpha+1/2)_k\,(-\gamma-n/2+1/2)_k\,k!}\,Q_{n-2k}(x;\gamma,\gamma,N)\;,$$

$$Q_n(x;\alpha,\beta,N) = \sum_{m=0}^{n} \frac{(\beta-\delta)_n\,(-1)^n}{(2+\alpha+\delta)_n}\cdot$$

$$\frac{(\alpha+\delta+1+2m)}{(\alpha+\delta+1)}\,\frac{(-n)_m\,(1+\alpha+\delta)_m\,(n+1+\alpha+\beta)_m\,(-1)^m}{(\alpha+2+n+\delta)_m\,(1-\beta+\delta-n)_m\,m!}\,Q_m(x;\alpha,\delta,N)$$

(compare [8], (4.1), (4.5)),

$$Q_n(x;\alpha,\beta,N) = \sum_{m=0}^{n} \frac{(\alpha-\gamma)_n\,(\beta+1)_n}{(\alpha+1)_n\,(2+\beta+\gamma)_n}\cdot$$

$$\frac{(\beta+\gamma+1+2m)}{(\beta+\gamma+1)}\,\frac{(-n)_m\,(1+\beta+\gamma)_m\,(\gamma+1)_m\,(n+1+\alpha+\beta)_m}{(\beta+1)_m\,(\beta+\gamma+n+2)_m\,(\gamma-\alpha-n+1)_m\,m!}\,Q_m(x;\gamma,\beta,N)$$

(compare [8], (4.1), (4.5)),

$$m_n^{(\gamma,\mu)}(x)=\sum_{m=0}^{n}\frac{(\gamma-\delta)_n\;(-n)_m}{(\delta-n+1-\gamma)_m\,m!}\;m_m^{(\delta,\mu)}(x)$$

(compare [8], (5.5)),

$$\widetilde{m}_n^{(\gamma,\mu)}(x)=\sum_{m=0}^{n}\left(\frac{\mu}{\mu-1}\right)^{n-m}\;\frac{(\gamma-\delta)_n\;(-n)_m}{(\delta-n+1-\gamma)_m\,m!}\;\widetilde{m}_m^{(\delta,\mu)}(x)\;,$$

$$m_n^{(\gamma,\mu)}(x)=\sum_{m=0}^{n}\left(\frac{\nu-\mu}{\mu\,(\nu-1)}\right)^{n}(\gamma)_n\;\frac{(-n)_m}{(\gamma)_m\,m!}\left(-\frac{\nu\,(\mu-1)}{\nu-\mu}\right)^{m}\;m_m^{(\gamma,\nu)}(x)$$

(compare [8], (5.4)),

$$\widetilde{m}_n^{(\gamma,\mu)}(x)=\sum_{m=0}^{n}\left(\frac{\nu-\mu}{(\mu-1)\,(\nu-1)}\right)^{n-m}\;\frac{(\gamma)_n\,(-n)_m\,(-1)^m}{(\gamma)_m\,m!}\;\widetilde{m}_m^{(\gamma,\nu)}(x)\;,$$

$$k_n^{(p)}(x,N)=\sum_{m=0}^{n}(p-q)^{n-m}\;\frac{(-N)_n\,(-n)_m\,(-1)^m}{n!\,(-N)_m}\;k_m^{(q)}(x,N)$$

(compare [8], (5.11)),

$$\widetilde{k}_n^{(p)}(x,N)=\sum_{m=0}^{n}(p-q)^{n-m}\;\frac{(-N)_n\,(-n)_m\,(-1)^m}{(-N)_m\,m!}\;\widetilde{k}_m^{(q)}(x,N)\;,$$

$$k_n^{(p)}(x,N)=\sum_{m=0}^{n}\frac{p^{n-m}\,(M-N)_n\,(-n)_m}{n!\,(N-M-n+1)_m}\;k_m^{(p)}(x,M)$$

(compare [8], (5.12)),

$$\widetilde{k}_n^{(p)}(x,N)=\sum_{m=0}^{n}\frac{p^{n-m}\,(M-N)_n\,(-n)_m}{(N-M-n+1)_m\,m!}\;\widetilde{k}_m^{(p)}(x,M)\;,$$

$$c_n^{(\mu)}(x)=\sum_{m=0}^{n}(-1)^n\;\frac{\nu^m}{\mu^n}\;(\nu-\mu)^{n-m}\;\frac{(-n)_m}{m!}\;c_m^{(\nu)}(x)$$

(compare [8], (5.16)),

$$\widetilde{c}_n^{(\mu)}(x)=\sum_{m=0}^{n}(-1)^m\;(\nu-\mu)^{n-m}\;\frac{(-n)_m}{m!}\;\widetilde{c}_m^{(\nu)}(x)\;,$$

$$K_n^{(\alpha,\beta)}(x)=\left(\frac{\beta-\delta}{\alpha}\right)_n\;\alpha^{n-m}\;\sum_{m=0}^{n}\frac{(-n)_m\;\left(\frac{\alpha(1-n)-\beta+\delta}{\alpha}\right)_m}{m!}\;K_m^{(\alpha,\delta)}(x)\;,$$

$$\widetilde{K}_n^{(\alpha,\beta)}(x)=\left(\frac{\beta-\delta}{\alpha}\right)_n\;\sum_{m=0}^{n}\frac{(-n)_m\;\left(\frac{\alpha(1-n)-\beta+\delta}{\alpha}\right)_m}{m!}\;\widetilde{K}_m^{(\alpha,\delta)}(x)\;.$$

*Proof:* Combining the representations

$$P_n(x) = \sum_{j\in\mathbb{Z}} A_j(n)\, x^{\underline{j}} \qquad \text{and} \qquad x^{\underline{j}} = \sum_{m\in\mathbb{Z}} B_m(j)\, Q_m(x)\ ,$$

and using Zeilberger's algorithm, the method of Corollary 2 yields the results.
The connection relations for the polynomials $K_n^{(\alpha,\beta)}(x)$ cannot be obtained by this method. Here Theorem 2 leads straightforwardly to the result. □

Although besides (56)–(58) the connection results were essentially known ([8], [3]), our development gives a unified treatment of them and makes new results like (56)–(58) easily accessible.
Note that some of the representations are rather complicated. We suggest the idea to use the notation ${}_pf_q$ for the summand of ${}_pF_q$, i.e.

$${}_pF_q\left(\begin{array}{c}\text{upper}\\ \text{lower}\end{array}\middle|\, x\right) = \sum_{k=0}^{\infty} {}_pf_q\left(\begin{array}{c}\text{upper}\\ \text{lower}\end{array}\middle|\, x;k\right) .$$

With this notation, (55) could be written in the standardized hypergeometric notation

$$h_n^{(\alpha,\beta)}(x,N) = {}_3f_1\left(\begin{array}{c}\beta-\delta, 1-N, \alpha+1\\ 2+\alpha+\delta\end{array}\middle|\, 1;n\right)\cdot$$

$$\cdot\sum_{m=0}^{n} {}_4f_5\left(\begin{array}{c}-n, 1+\alpha+\delta, n+1+\alpha+\beta, \alpha/2+\delta/2+3/2, 1\\ 1-N, \alpha+1, \alpha/2+\delta/2+1/2, \alpha+2+n+\delta, -n-\beta+\delta+1\end{array}\middle|\, 1;m\right) h_m^{(\alpha,\delta)}(x,N)\ .$$

Finally, we deduce the parameter derivatives for the classical discrete orthogonal polynomials.

**Corollary 5** The following representations for the parameter derivatives of the classical discrete orthogonal polynomials are valid:

$$\begin{aligned}\frac{\partial}{\partial\alpha} h_n^{(\alpha,\beta)}(x,N) &= \sum_{m=0}^{n-1} \frac{1}{\alpha+\beta+m+n+1}\cdot\Big(h_n^{(\alpha,\beta)}(x,N) + \\ &\frac{(-1)^{n-m}\,(\alpha+\beta+1+2m)\,(1-N+m)_{n-m}\,(\beta+1+m)_{n-m}}{(n-m)\,(\alpha+\beta+1+m)_{n-m}}\, h_m^{(\alpha,\beta)}(x,N)\Big)\ ,\end{aligned}$$

$$\frac{\partial}{\partial\alpha}\tilde{h}_n^{(\alpha,\beta)}(x,N) = \sum_{m=0}^{n-1} \frac{(-1)^{n-m}\,(\alpha+\beta+1+2m)}{(\alpha+\beta+m+n+1)\,(n-m)}\,\frac{(1-N+m)_{n-m}\,(\beta+1+m)_{n-m}\,n!}{(\alpha+\beta+1+2m)_{2n-2m}\,m!}\,\tilde{h}_m^{(\alpha,\beta)}(x,N)\,,$$

$$\begin{aligned}\frac{\partial}{\partial\alpha} Q_n(x;\alpha,\beta,N) &= \sum_{m=0}^{n-1}\left(\frac{1}{\alpha+\beta+m+n+1} - \frac{1}{\alpha+m+1}\right)\cdot\Big(Q_n(x;\alpha,\beta,N) + \\ &\frac{(\alpha+\beta+1+2m)\,(\beta+1+m)_{n-m}\,n!}{(n-m)\,(\alpha+1+m)_{n-m}\,(\alpha+\beta+1+m)_{n-m}\,m!}\, Q_m(x;\alpha,\beta,N)\Big)\ ,\end{aligned}$$

$$\begin{aligned}\frac{\partial}{\partial\beta} h_n^{(\alpha,\beta)}(x,N) &= \sum_{m=0}^{n-1} \frac{1}{\alpha+\beta+m+n+1}\cdot\Big(h_n^{(\alpha,\beta)}(x,N) + \\ &\frac{\alpha+\beta+1+2m}{n-m}\,\frac{(1-N+m)_{n-m}\,(\alpha+1+m)_{n-m}}{(\alpha+\beta+1+m)_{n-m}}\, h_m^{(\alpha,\beta)}(x,N)\Big)\ ,\end{aligned}$$

$$\frac{\partial}{\partial \beta}\widetilde{h}_n^{(\alpha,\beta)}(x,N)=\sum_{m=0}^{n-1}\frac{\alpha+\beta+1+2m}{(\alpha+\beta+m+n+1)\,(n-m)}\,\frac{(1-N+m)_{n-m}\,(\alpha+1+m)_{n-m}\,n!}{(\alpha+\beta+1+2m)_{2n-2m}\,m!}\,\widetilde{h}_m^{(\alpha,\beta)}(x,N)\,,$$

$$\begin{aligned}\frac{\partial}{\partial \beta}Q_n(x;\alpha,\beta,N) \quad = \quad & \sum_{m=0}^{n-1}\frac{1}{\alpha+\beta+m+n+1}\cdot\Big(Q_n(x;\alpha,\beta,N)+ \\ & \frac{(-1)^{n-m}\,(\alpha+\beta+1+2m)}{n-m}\,\frac{n!}{(\alpha+\beta+1+m)_{n-m}\,m!}\,Q_m(x;\alpha,\beta,N)\Big)\ ,\end{aligned}$$

$$\frac{\partial}{\partial \mu}m_n^{(\gamma,\mu)}(x) = \frac{n(-1+\gamma+n)}{(1-\mu)\mu}\,m_{n-1}^{(\gamma,\mu)}(x) - \frac{n}{(1-\mu)\mu}\,m_n^{(\gamma,\mu)}(x)$$

(see e.g. [12], Theorem 9),

$$\frac{\partial}{\partial \mu}\widetilde{m}_n^{(\gamma,\mu)}(x) = \frac{n(1-\gamma-n)}{(1-\mu)^2}\,\widetilde{m}_{n-1}^{(\gamma,\mu)}(x)\ ,$$

$$\frac{\partial}{\partial \gamma}m_n^{(\gamma,\mu)}(x) = \sum_{m=0}^{n-1}\frac{n!}{m!\,(n-m)}\,m_m^{(\gamma,\mu)}(x)$$

(see [12], Theorem 10),

$$\frac{\partial}{\partial \gamma}\widetilde{m}_n^{(\gamma,\mu)}(x) = \sum_{m=0}^{n-1}\left(\frac{\mu}{\mu-1}\right)^{n-m}\frac{n!}{m!\,(n-m)}\,\widetilde{m}_m^{(\gamma,\mu)}(x)\ ,$$

$$\frac{\partial}{\partial p}k_n^{(p)}(x,N) = (-1+n-N)\,k_{n-1}^{(p)}(x,N)$$

(see e.g. [12], Theorem 9),

$$\frac{\partial}{\partial p}\widetilde{k}_n^{(p)}(x,N) = n\,(-1+n-N)\,\widetilde{k}_{n-1}^{(p)}(x,N)\ ,$$

$$\frac{\partial}{\partial \mu}c_n^{(\mu)}(x) = \frac{n}{\mu}\,c_{n-1}^{(\mu)}(x) - \frac{n}{\mu}\,c_n^{(\mu)}(x)$$

(see e.g. [12], Theorem 9),

$$\frac{\partial}{\partial \mu}\widetilde{c}_n^{(\mu)}(x) = -n\,\widetilde{c}_{n-1}^{(\mu)}(x)\ ,$$

$$\frac{\partial}{\partial \beta}K_n^{(\alpha,\beta)}(x) = \sum_{m=0}^{n-1}\frac{\alpha^{n-m-1}\,n!}{(n-m)\,m!}\,K_m^{(\alpha,\beta)}(x)\ ,$$

$$\frac{\partial}{\partial \beta}\widetilde{K}_n^{(\alpha,\beta)}(x) = \sum_{m=0}^{n-1}\frac{n!}{\alpha\,(n-m)\,m!}\,\widetilde{K}_m^{(\alpha,\beta)}(x)\ .$$

*Proof:* If the derivative is taken with respect to a variable occurring as an *argument* rather than as a *parameter* in the hypergeometric representation, its representation can be easily obtained from the derivative rule of the generalized hypergeometric function, and the chain rule. In those cases, the representations need at most two neighboring polynomials.
The other cases can be handled similarly to Corollary 3. □

# 7 Conclusion

Here, we want to recall the algorithms to convert between the different types of representations:

1. **Hypergeometric Representation → Recurrence Equation**: Zeilberger's algorithm
2. **Hypergeometric Representation → Difference/Differential Equation**: Zeilberger's/Almkvist-Zeilberger's algorithm
3. **Difference/Differential Equation → Recurrence Equation**: Theorem 1
4. **Difference/Differential Equation → Hypergeometric Representation**: method of § 3 and § 5
5. **Recurrence Equation → Difference/Differential Equation**: Algorithms 1 and 2 in [14]
6. **Recurrence Equation → Hypergeometric Representation**: combination of methods 5 and 4